\documentclass[10pt,reqno]{amsart}

\usepackage{fullpage}
\usepackage{amsmath,amsthm,amssymb,thmtools,tikz-cd,comment,hyperref,cleveref,enumitem,breqn,float,caption,subfig,enumitem,ytableau,tikz,cases,float,kotex,multirow}
\usepackage{xcolor}
\usetikzlibrary{positioning, arrows.meta, calc}

\usetikzlibrary{decorations.pathreplacing,
	calligraphy,
	matrix}

\allowdisplaybreaks
\theoremstyle{plain}
\newtheorem{theorem}{Theorem}[section]
\newtheorem{lemma}[theorem]{Lemma}
\newtheorem{proposition}[theorem]{Proposition}
\newtheorem{corollary}[theorem]{Corollary}
\newtheorem{thmalphabetintro}{Theorem}

\newtheorem{thmalphabetmaintext}{Theorem}

\newtheorem{coralphabetintro}[thmalphabetintro]{Corollary}

\newtheorem{coralphabetmaintext}[thmalphabetmaintext]{Corollary}

\theoremstyle{definition}
\newtheorem{definition}[theorem]{Definition}
\newtheorem{example}[theorem]{Example}

\newtheorem{remark}[theorem]{Remark}
\newtheorem{conjecture}[theorem]{Conjecture}
\newtheorem{question}[theorem]{Question}
\newtheorem{notation}[theorem]{Notation}

\newtheorem{open problem}[theorem]{Open Problem}
\newtheorem{remark and notation}[theorem]{Remark and Notation}
\newtheorem{remark and definition}[theorem]{Remark and Definition}
\newtheorem{definition and notation}[theorem]{Definition and Notation}
\newtheorem{notation and convention}[theorem]{Notation and convention}
\newtheorem{convention and notation}[theorem]{Convention and notation}

\def \p {\mathbb{P}}

\def \o {\mathcal{O}}
\def \i {\mathcal{I}}

\def \l {\mathcal{L}}

\def \N {\textnormal{\bf{N}}}
\def \H {\textnormal{H}}

\def \proj {\operatorname{Proj}}
\def \codim {\operatorname{codim}}

\def \tor {\operatorname{Tor}}

\def \gon {\operatorname{gon}}
\def \vmd[#1]{\text{VMD$^{#1}$}}
\def \dpv[#1]{\text{dP$^{#1}$}}
\def \acm[#1,#2]{\text{ACMD$^{#1}_{#2}$}}
\def \propp[#1,#2]{\text{$P_{{#1},{#2}}$}}
\def \propa[#1,#2]{\text{$A_{{#1},{#2}}$}}

\def\p{{\mathbb P}}

\title[Hierarchical structure of quadratic strand]{Hierarchical structure of graded Betti numbers in the quadratic strand}
\author{Jong In Han}
\address{Jong In Han, School of Mathematics, Korea Institute for Advanced Study (KIAS), 85 Hoegi-ro, Dongdaemun-gu, Seoul, 02455, Republic of Korea}
\email{jihan09@kias.re.kr}

\author{Sijong Kwak}
\address{Sijong Kwak, Department of Mathematical Sciences, Korea Advanced Institute of Science and Technology (KAIST), 291, Daehak-ro, Yuseong-gu, Daejeon, Republic of Korea}
\email{sjkwak@kaist.ac.kr}

\author{Wanseok Lee}
\address{Wanseok Lee, Department of Applied Mathematics, Pukyong National University, 45, Yongso-ro, Nam-gu, Busan, Republic of Korea}
\email{wslee@pknu.ac.kr}

\thanks{}

\begin{document}
	\begin{abstract}
		The classical results, initiated by Castelnuovo and Fano and later refined by Eisenbud and Harris, provide several upper bounds on the number of quadrics defining a nondegenerate projective variety. Recently, it has been revealed that these bounds extend naturally to certain linear syzygies, suggesting the presence of a hierarchical structure governing the quadratic strand of graded Betti numbers. 
		
		In this article, we establish such a hierarchy in full generality. We first prove sharp upper bounds for $\beta_{p,1}(X)$ depending on the degree of a projective variety $X$, extending the classical quadratic bounds to all linear syzygies and identifying the extremal varieties in each range.
		We then introduce geometric conditions that describe how containment of $X$ in low-degree varieties influences syzygies, and we show that these conditions stratify the quadratic strand into a finite sequence of hierarchies. This leads to a complete description of all possible extremal behavior. We also prove a generalized $K_{p,1}$-theorem, demonstrating that the vanishing of $\beta_{p,1}(X)$ detects containment in a variety of minimal degree at each hierarchy. 
	\end{abstract}
	\maketitle
	
	\section{Introduction}
	Let $X\subseteq\p^r$ be a nondegenerate projective variety of dimension $n\geq 1$, codimension $e$, and degree $d$ over an algebraically closed field $\Bbbk$ of characteristic zero. 
	Let $S = \Bbbk[x_0, \ldots, x_r]$ be the homogeneous coordinate ring of $\p^r$. Denote the homogeneous ideal of $X$ by $I_X$ and the homogeneous coordinate ring by $S_X$. 
	The \textit{$(i,j)$-th graded Betti number of $X$} is defined as $$\beta_{i,j}(X)=\dim_\Bbbk\tor_i(S_X,\Bbbk)_{i+j}.$$ 
	We say $X$ satisfies the property $\N_{2,p}$ if $\beta_{i,j}(X)=0$ for all $i\leq p$ and $j\geq 2$.
	The Betti number $\beta_{i,j}(X)$ is said to lie in the quadratic strand if $j=1$. 
	If $d=e+1$, then we call $X$ a \textit{variety of minimal degree} (VMD).
	If $d=e+2$, then we call $X$ a \textit{variety of almost minimal degree}.
	If, in addition, $X$ is arithmetically Cohen--Macaulay (ACM), then it is called a \textit{del Pezzo variety}.
	For convenience, we use the following notation throughout the paper.
	\begin{notation}
		We denote a variety of minimal degree (resp. a del Pezzo variety) of dimension $n$ by \vmd[n] (resp. \dpv[n]).
		We also denote by \acm[k,m] an integral ACM divisor $Y$ of dimension $k$ with $\deg(Y)=\codim(Y)+m$ contained in \vmd[k+1].
	\end{notation}
	
	It is a classical problem to study the number $\beta_{1,1}(X)$ of quadratic forms defining $X$ and the syzygies among them. By G. Castelnuovo \cite{zbMATH02692308} and G. Fano \cite{zbMATH02683130}, the following classical results hold.
	
	\begin{itemize}
		\item It always holds that $\beta_{1,1}(X)\leq \binom{e+1}{2}$. The equality holds if and only if $X$ is a variety of minimal degree, i.e., $X$ is a $\vmd[n]$.
		\item If $X$ is not a variety of minimal degree, then $\beta_{1,1}(X)\leq \binom{e+1}{2}-1$. The equality holds if and only if $X$ is a $\dpv[n]$ (when $e\geq 2$).
	\end{itemize}
	
	The Castelnuovo theory developed in \cite{zbMATH02692308, zbMATH02683130, MR685427} exhibits a similar staircase phenomenon in higher degrees. 
	Subsequent work by Petrakiev~\cite{MR2414947} and Park~\cite{MR4145215} provided a deeper understanding of the behavior of $ \beta_{1,1}(X) $ beyond the minimal and del Pezzo cases. Specifically, the following holds.
	
	\begin{itemize}
		\item If $d\geq 2e+1$, then $\beta_{1,1}(X)\leq \binom{e}{2}$. 
		
		\noindent
		If $d\geq 2e+3$, then the equality holds if and only if $X$ is contained in a \vmd[n+1] ($e\geq 2$).
		\item Hence, if $d\geq 2e+3$ and $X$ is not contained in a \vmd[n+1], then $\beta_{1,1}(X)\leq\binom{e}{2}-1$ ($e\geq 2$).
		
		\noindent
		If $d\geq 2e+5$, then the equality holds if and only if $X$ is contained in a \dpv[n+1] ($e\geq 3$).
		
		\noindent
		\item For $k \in \{2,3\}$, suppose that $e \geq	4$ if $k=2$, and $e=6$ if $k=3$. If $d \geq 2e+2k+3$, then the equality $\beta_{1,1}(X)=\binom{e}{2}-k$ holds if and only if there exists a unique $(n+1)$-dimensional variety $Y\subseteq\mathbb{P}^{r}$ such that $X\subseteq Y$ and $\beta_{1,1}(Y)=\binom{e}{2}-k$.
	\end{itemize}
	
	In addition, sharp upper bounds for the Betti numbers in the quadratic strand have been established in \cite{MR3302629, MR4645613}.
	\begin{itemize}
		\item It always holds that $\beta_{p,1}(X)\leq p\binom{e+1}{p+1}$ for all $p\geq 1$, and the equality holds for all $p$ if and only if $X$ is a \vmd[n].
		\item If $X$ is not a variety of minimal degree, then $\beta_{p,1}(X)\leq p\binom{e+1}{p+1}-\binom{e}{p-1}$ $(1\leq p\leq e)$.
		
		\noindent
		The equality holds for all $p$ if and only if $X$ is a \dpv[n] ($e\geq 2$).
		\item If $e\geq 2$ and X is neither a \vmd[n] nor a \dpv[n], then $\beta_{p,1}(X)\leq p\binom{e+1}{p+1}-2\binom{e}{p-1}$ $(1\leq p\leq e-1)$.
		
		\noindent
		The equality holds for all $p$ if and only if $X$ is either a variety of almost minimal degree with arithmetic depth $n$ and the Green--Lazarsfeld index zero, or an ACM variety of degree $d=e+3$ ($e\geq 3$).
	\end{itemize}
	
	The following table provides a concise summary of the known results.
	
	\begin{table}[H]
		\resizebox{\linewidth}{!}{%
			\begin{tabular}{|c||c|c||cccc|c|}
				\hline
				Condition                                                                      & $\beta_{1,1}$          & Extremal case        & \multicolumn{1}{c|}{\quad$\beta_{2,1}$\quad\quad} & \multicolumn{1}{c|}{\quad$\beta_{3,1}$\quad\quad} & \multicolumn{1}{c|}{\quad\quad\quad\quad\quad$\cdots$\quad\quad\quad\quad\quad\quad} & \quad$\beta_{e,1}$\quad\quad & Extremal case \\ \hline\hline
				$d\geq e+1$                                                                    & $\leq\binom{e+1}{2}$\quad\quad   & \vmd[n]              & \multicolumn{4}{c|}{$\leq p\binom{e+1}{p+1}$\quad\quad\quad\quad}                                                                           & $\vmd[n]$       \\ \hline
				$d\geq e+2$                                                                    & $\leq\binom{e+1}{2}-1$ & \dpv[n]              & \multicolumn{4}{c|}{$\leq p\binom{e+1}{p+1}-\binom{e}{p-1}$}                                                            & $\dpv[n]$       \\ \hline
				$d\geq e+3$                                                                    & $\leq\binom{e+1}{2}-2$ & ACM, $d=e+3$          & \multicolumn{4}{c|}{$\leq p\binom{e+1}{p+1}-2\binom{e}{p-1}$}                                                           & ACM, d=e+3   \\ \hline
				$d\geq e+4$                                                                    & $\leq\binom{e+1}{2}-3$ & ACM, $d=e+4$          & \multicolumn{4}{c|}{unknown}                                                           & -   \\ \hline
				$d\geq e+5$                                                                    & $\leq\binom{e+1}{2}-4$ & ACM, $d=e+5$          & \multicolumn{4}{c|}{unknown}                                                           & -   \\ \hline
				$\vdots$                                                                    & $\vdots$ & $\vdots$          & \multicolumn{4}{c|}{$\vdots$}                                                           & $\vdots$   \\ \hline
				$d\geq 2e$                                                                    & $\leq\binom{e+1}{2}-(e-1)$ & ACM, $d=2e$          & \multicolumn{4}{c|}{unknown}                                                           & -   \\ \hline
				$d\geq 2e+1$                                                                   & $\leq\binom{e}{2}$     & not classified              & \multicolumn{4}{c|}{unknown}                                                                                            & -             \\ \hline
				$d\geq 2e+3$                                                                   & $\leq\binom{e}{2}$     & $\subseteq$\vmd[n+1] & \multicolumn{4}{c|}{unknown}                                                                                            & -             \\ \hline
				\begin{tabular}[c]{@{}c@{}}$d\geq 2e+3$\\ $X\nsubseteq\vmd[n+1]$\end{tabular} & $\leq\binom{e}{2}-1$   & not classified              & \multicolumn{4}{c|}{unknown}                                                                                            & -             \\ \hline
				\begin{tabular}[c]{@{}c@{}}$d\geq 2e+5$\\ $X\nsubseteq\vmd[n+1]$\end{tabular} & $\leq\binom{e}{2}-1$   & $\subseteq$\dpv[n+1] & \multicolumn{4}{c|}{unknown}                                                                                            & -             \\ \hline
			\end{tabular}%
		}
		\vspace{3pt}
		\caption{Well-known upper bounds of $\beta_{p,1}$}
		\label{table:previous bounds}
	\end{table}
	
	The first three rows in the table exhibit a possibility that there might exist a hierarchical structure for graded Betti numbers $\beta_{p,1}$ in the quadratic strand.
	Our main purpose is to show the existance of such a hierarchical structure.
	First, we establish the first hierarchy.
	\begin{thmalphabetintro}[The first hierarchy]\label{hierarchy_upper_bound_intro}
		Let $X\subseteq\p^r$ be a nondegenerate projective variety of dimension $n$, codimension $e\geq 2$, and degree $d\geq e+1+m$ for an integer $0\leq m\leq e$. Then we have
		\begin{numcases}{\beta_{p,1}(X)\leq}
			p\binom{e+1}{p+1}-m\binom{e}{p-1} &for $1\leq p\leq e+1-m$\\
			p\binom{e}{p+1} &for $e+1-m\leq p$.
		\end{numcases}
		When $2\leq m\leq e-1$, the equality holds for all $p$ if and only if $X$ is an ACM divisor of degree $e+m+1$ in a variety of minimal degree, i.e., $X$ is an \acm[n,m+1].
	\end{thmalphabetintro}
	
	This reveals the existence of a hierarchy in the Betti numbers with respect to the degree. 
	Note that the inequality (2) essentially appears from $m \ge 3$ (see \Cref{table:first hierarchy}), a phenomenon that had not been observed in the previous works~\cite{MR3302629, MR4645613}. When $p=e+1-m$, the bounds~(1) and~(2) coincide.
	
	Furthermore, the following corollary extends the classical Castelnuovo bound on the number of quadrics that asserts $\beta_{1,1}\leq \binom{e}{2}$ when $d\geq 2e+1$.
	\begin{coralphabetintro}\label{upper_bound_second_floor_intro}
		With the same notation as above, when $d\geq 2e+1$, we have \[\beta_{p,1}(X)\leq p\binom{e}{p+1}\] for all $p\geq 1$.
		The extremal varieties are exactly the ones contained in a \vmd[n+1].
	\end{coralphabetintro}
	
	The following table summarizes \Cref{hierarchy_upper_bound_intro} and \Cref{upper_bound_second_floor_intro}.
	In the table, $\pi$ represents the sectional genus of the extremal case.
	
	\begin{table}[H]
		\resizebox{\linewidth}{!}{%
			\begin{tabular}{|c||ccccccccc|c|c|c|}
				\hline
				Condition    & \multicolumn{1}{c|}{$\beta_{1,1}$} & \multicolumn{1}{c|}{$\beta_{2,1}$} & \multicolumn{1}{c|}{$\beta_{3,1}$} & \multicolumn{1}{c|}{$\cdots$} & \multicolumn{1}{c|}{$\beta_{e-4,1}$} & \multicolumn{1}{c|}{$\beta_{e-3,1}$} & \multicolumn{1}{c|}{$\beta_{e-2,1}$} & \multicolumn{1}{c|}{$\beta_{e-1,1}$}        & $\beta_{e,1}$      & Extremal case          & $\pi$    & Reference                                                  \\ \hline\hline
				$d\geq e+1$  & \multicolumn{9}{c|}{$\leq p\binom{e+1}{p+1}$\quad\quad\quad\quad}                                                                                                                                                                                                                                                                                        & \vmd[n]                & 0        & \cite{MR3302629}                                           \\ \hline
				$d\geq e+2$  & \multicolumn{8}{c|}{$\leq p\binom{e+1}{p+1}-\binom{e}{p-1}$}                                                                                                                                                                                                                                                    & \multirow{8}{*}{0} & \dpv[n]                & 1        & \cite{MR3302629}                                           \\ \cline{1-9} \cline{11-13} 
				$d\geq e+3$  & \multicolumn{8}{c|}{$\leq p\binom{e+1}{p+1}-\textcolor{blue}{2}\binom{e}{p-1}$\quad\quad\quad}                                                                                                                                                                                                                                                   &                    & \acm[n,3]            & 2        & \cite{MR4645613}                                           \\ \cline{1-9} \cline{11-13} 
				$d\geq e+4$  & \multicolumn{7}{c|}{$\leq p\binom{e+1}{p+1}-\textcolor{blue}{3}\binom{e}{p-1}$}                                                                                                                                                                                                     & \multicolumn{1}{c|}{$\leq p\binom{e}{p+1}$} &                    & \acm[n,4]            & 3        & \multirow{6}{*}{\Cref{hierarchy_upper_bound_intro}}        \\ \cline{1-9} \cline{11-12}
				$d\geq e+5$  & \multicolumn{6}{c|}{$\leq p\binom{e+1}{p+1}-\textcolor{blue}{4}\binom{e}{p-1}$}                                                                                                                                                              & \multicolumn{2}{c|}{$\leq p\binom{e}{p+1}$}                                        &                    & \acm[n,5]            & 4        &                                                            \\ \cline{1-9} \cline{11-12}
				$d\geq e+6$  & \multicolumn{5}{c|}{$\leq p\binom{e+1}{p+1}-\textcolor{blue}{5}\binom{e}{p-1}$}                                                                                                                       & \multicolumn{3}{c|}{$\leq p\binom{e}{p+1}$}                                                                               &                    & \acm[n,6]            & 5        &                                                            \\ \cline{1-9} \cline{11-12}
				$\vdots$     & \multicolumn{8}{c|}{$\vdots$}                                                                                                                                                                                                                                                                                   &                    & $\vdots$               & $\vdots$ &                                                            \\ \cline{1-9} \cline{11-12}
				$d\geq 2e$   & \multicolumn{2}{c|}{$\leq p\binom{e+1}{p+1}-\textcolor{blue}{(e-1)}\binom{e}{p-1}$}       & \multicolumn{6}{c|}{$\leq p\binom{e}{p+1}$}                                                                                                                                                                                           &                    & \acm[n,e]            & $e-1$    &                                                            \\ \cline{1-9} \cline{11-13} 
				$d\geq 2e+1$ & \multicolumn{8}{c|}{$\leq p\binom{e}{p+1}$}                                                                                                                                                                                                                                                                     &                    & $\subseteq\vmd[n+1]$ & -        & \multicolumn{1}{l|}{\Cref{upper_bound_second_floor_intro}} \\ \hline
			\end{tabular}%
		}
		\vspace{3pt}
		\caption{The first hierarchy}
		\label{table:first hierarchy}
	\end{table}
	
	The following theorem provides a quadratic-strand analogue of the Castelnuovo theory on the number of quadrics, asserting that $\beta_{1,1}(X)\leq \binom{e}{2}-1$ when $d\geq 2e+3$ and $X$ is not contained in a $\vmd[n+1]$.
	
	\begin{thmalphabetintro}[The next row]\label{upper_bound_second_floor_next_intro}
		Let $X\subseteq\p^r$	 be a nondegenerate projective variety of dimension $n$, codimension $e\geq 2$, and degree $d$ that is not contained in a \vmd[n+1]. If $d\gg 0$, then 
		\[\beta_{p,1}(X)\leq p\binom{e}{p+1}-\binom{e-1}{p-1}\]
		for all $p\geq 1$. Furthermore, the extremal varieties are the ones contained in a \dpv[n+1] when $e\geq 3$.
	\end{thmalphabetintro}
	
	\Cref{upper_bound_second_floor_next_intro} adds an additional row to \Cref{table:first hierarchy} which gives a hint that there might exist the next hierarchy.
	
	\begin{table}[H]
		\resizebox{\linewidth}{!}{%
			\begin{tabular}{|c||ccccccc|cc|c|c|}
				\hline
				Condition                                                                  & \multicolumn{1}{c|}{$\beta_{1,1}$} & \multicolumn{1}{c|}{$\beta_{2,1}$} & \multicolumn{1}{c|}{$\beta_{3,1}$} & \multicolumn{1}{c|}{$\cdots$} & \multicolumn{1}{c|}{$\beta_{e-4,1}$} & \multicolumn{1}{c|}{$\beta_{e-3,1}$} & $\beta_{e-2,1}$ & \multicolumn{1}{c|}{$\beta_{e-1,1}$} & $\beta_{e,1}$ & Extremal case         & Reference                                  \\ \hline\hline
				\begin{tabular}[c]{@{}c@{}}$d\gg 0$\\ $X\nsubseteq \vmd[n+1]$\end{tabular} & \multicolumn{7}{c|}{$\leq p\binom{e}{p+1}-\binom{e-1}{p-1}$}                                                                                                                                                                                      & \multicolumn{2}{c|}{0}                               & $\subseteq \dpv[n+1]$ & \Cref{upper_bound_second_floor_next_intro} \\ \hline
			\end{tabular}%
		}
		\vspace{3pt}
		\caption{The next upper bounds with the extremal case}
		\label{table:to the next hierarchy}
	\end{table}
	
	One might expect that the same pattern in the first hierarchy continues again.
	In other words, one may expect if $d\gg 0$ and $X\nsubseteq\vmd[n+1]$ and $X\nsubseteq\dpv[n+1]$, then $\beta_{p,1}\leq p\binom{e}{p+1}-2\binom{e-1}{p-1}$.
	Indeed, this holds true and the complete form of the second hierarchy exists. 
	This phenomenon does not stop here, and there exist further hierarchies (see \Cref{hierarchy_theorem_intro} and \Cref{table:complete hierarchies}).
	To achieve these, we revisit the Castelnuovo theory in more detail.
	First, we introduce the geometric conditions $\propa[k,m]$ and $\propp[k,m]$ as follows.
	\begin{definition}
		Let $X\subseteq\p^r$ be a nondegenerate algebraic set of dimension $n$, codimension $e$, and degree $d$.
		For nonnegative integers $k,m$, we use the following definition:
		\begin{itemize}
			\item[\textnormal{(i)}] $X$ satisfies $\propa[k,m]$ if $X$ is not contained in any variety $Y\subseteq\p^r$ of dimension $n+k$ with $\deg Y\leq \codim Y+m$ (i.e., whenever $X\subseteq Y$ for a variety $Y$ of dimension $n+k$, we have $\deg Y\geq \codim Y+m+1$);
			\item[\textnormal{(ii)}] $X$ satisfies $\propp[k,m]$ if $X$ satisfies $\propa[k,m]$ and $\propa[i,d_i]$ for all $0\leq i<k$ where $d_i:=2^{e-i}-(e-i)$.
		\end{itemize}
		Also, for convenience, when $k<0$, we define $\propp[k,m]$ to be the vacuous condition, i.e., it is always true for any $X$. However, we almost always assume $k$ and $m$ are nonnegative. The only case we consider that $k<0$ is in \Cref{Gen_Kp1_intro}.
	\end{definition}
	
	For example,
	\begin{enumerate}
		\item $X$ satisfies $\propa[0,m]$ if $d\geq e+m+1$;
		\item $X$ satisfies $\propa[1,1]$ if it is not contained in a \vmd[n+1];
		\item $X$ satisfies $\propa[1,2]$ if it is not contained in an $(n+1)$-fold $Y$ with $\deg Y\leq \codim Y+2$;
		\item $X$ satisfies $\propa[2,1]$ if it is not contained in a \vmd[n+2];
		\item The properties $\propp[0,m]$ and $\propa[0,m]$ are the same.
		\item $X$ satisfies $\propp[1,1]$ if $d>2^e$ and $X$ is not contained in a \vmd[n+1];
		\item $X$ satisfies $\propp[1,2]$ if $d>2^e$ and $X$ is not contained in an $(n+1)$-fold $Y$ with $\deg Y\leq \codim Y+2$;
		\item \propp[k,m+1] implies \propp[k,m]. 
		\item \propp[k+1,0] implies \propp[k,m] for all $m\leq d_{k}$.
	\end{enumerate}
	
	On the other hand, together with Eisenbud, Harris~\cite{MR685427} developed a modern refinement of
	Castelnuovo’s classical theory by introducing a decreasing sequence $\pi_m(d,r)$ for $m=0,1,\ldots,r-1$.
	These numbers provide refined upper bounds for the arithmetic genus of curves in $\mathbb P^r$ under suitable containment constraints (e.g., not lying on low-degree surfaces). This led to the following conjecture, which may be viewed as 
	a Castelnuovo-type analogue for curves of nearly maximal genus.
	
	\begin{conjecture}[{\cite{MR685427}, also cf. \cite{MR2414947}}]\label{contraposition_conj_ver1}
		Let $ C \subseteq \mathbb{P}^r $ be a nondegenerate projective curve of degree $ d \ge 2r + 2m - 1 $. 
		If $ C $ satisfies $\propa[1,m]$, then $ g \le \pi_m(d,r) $.
	\end{conjecture}
	
	Petrakiev~\cite{MR2414947} approached this conjecture from a different perspective by translating the genus-bound problem for curves into a statement about the quadratic relations of finite sets of points in symmetric position. 
	The Eisenbud--Harris conjecture can be reformulated in terms of $\beta_{1,1}$ for finite sets of points, thus reducing the geometric complexity to a more algebraic problem in syzygy theory. This led to the following conjecture, which can be seen as an algebraic counterpart to Conjecture~\ref{contraposition_conj_ver1}.
	
	\begin{conjecture}[{\cite{MR2414947}}]\label{contraposition_conj_ver2}
		Let $ \Gamma \subseteq \mathbb{P}^{e} $ be a set of $ d \ge 2e + 2m + 1 $ points in symmetric position $(1 \le m \le e - 2)$. 
		If $ \Gamma $ satisfies $\propa[1,m]$, then 
		\begin{equation*}
			\beta_{1,1}(\Gamma) \neq \binom{e}{2} +1 - m .
		\end{equation*}
	\end{conjecture}
	
	This conjecture remains open. However, Eisenbud and Harris \cite{MR685427} proved that Conjecture~\ref{contraposition_conj_ver1} 
	holds with an additional condition $d > d_0'$ where
	\[
	d_0'=\left\{
	\begin{array}{ll}
		36r &\text{when }r\leq 6,\\
		288 &\text{when }r=7,\\
		2^{r+1} &\text{when }r\geq 8.
	\end{array}
	\right.
	\]
	Petrakiev’s reformulation played a crucial role in shifting the focus from geometric genus bounds to algebraic properties of syzygies, laying the groundwork for the hierarchical structures developed in this paper.
	
	\begin{theorem}[{\cite{MR685427}}]
		If $C$ satisfies $d > d_0'$ and $\propa[1,m]$, then \Cref{contraposition_conj_ver1} is true.
	\end{theorem}

	Along this line, we first show that \Cref{contraposition_conj_ver2} holds for projective varieties under the condition $\propp[1,m]$ (i.e., $d > d_0$ and $\propa[1,m]$), which follows from the case $k=1$ in the following \Cref{hierarchy_theorem_intro}. Moreover, we extend the upper bounds of $\beta_{1,1}$ to  those of higher linear syzygies $\beta_{p,1}$, and these exactly correspond to the second hierarchy. However, this phenomenon does not terminate here; rather, it unfolds into further layers of hierarchies.

	\begin{thmalphabetintro}[The complete hierarchy]\label{hierarchy_theorem_intro}
		Let $X\subseteq\p^r$ be a nondegenerate projective variety of dimension $n$, codimension $e\geq 2$, and degree $d$. For some $0\leq k\leq e-1$ and $0\leq m\leq e-k$, if X satisfies $\propp[k,m]$, then
		\begin{numcases}{\beta_{p,1}(X)\leq}
			p\binom{e+1-k}{p+1}-m\binom{e-k}{p-1} &for $1\leq p\leq e+1-m-k$ \nonumber \\
			p\binom{e-k}{p+1} &for $e+1-m-k\leq p$. \nonumber
		\end{numcases}
		Furthermore, the extremal cases are determined as follows except $(k,m)=(e-1,1)$.
		\begin{itemize}
			\item[\textnormal{(i)}] When $m=0$, the equalities hold for all $p\geq 1$ if and only if $X$ is contained in a \vmd[n+k].
			\item[\textnormal{(ii)}] When $m=1$, the equalities hold for all $p\geq 1$ if and only if $X$ is contained in a \dpv[n+k].
			\item[\textnormal{(iii)}] When $2\leq m\leq e-k-1$, the equalities hold for all $p\geq 1$ if and only if $X$ is contained in an \acm[n+k,m+1].
			\item[\textnormal{(iv)}] When $m=e-k$, the equalities hold for all $p\geq 1$ if and only if $X$ is contained in a \vmd[n+k+1].
		\end{itemize}
		When $(k,m)=(e-1,1)$, all varieties (satisfying \propp[e-1,1]) are vacuously extremal.
	\end{thmalphabetintro}
	
	Setting $k=0$ in \Cref{hierarchy_theorem_intro} recovers the bounds of \Cref{hierarchy_upper_bound_intro} and \Cref{upper_bound_second_floor_intro}.
	The special case of $k=1$ and $m=1$ in \Cref{hierarchy_theorem_intro} reduces to \Cref{upper_bound_second_floor_next_intro}.
	Furthermore, \Cref{hierarchy_theorem_intro} establishes a generalized version of the $K_{p,1}$-theorem (\Cref{Gen_Kp1_intro}).
	\Cref{hierarchy_theorem_intro} is summarized in the following table.
	Remind that \propp[k,m+1] implies \propp[k,m], and if $m\leq d_{k}$, then \propp[k+1,0] implies \propp[k,m].
	That is, if any condition in the table holds, the condition above it necessarily holds.

	\begin{table}[H]
		\resizebox{\linewidth}{!}{%
			\begin{tabular}{c|c||ccccccccc|c|c|}
				\cline{2-13}
				\multicolumn{1}{l|}{}                                       & Condition                                                                             & \multicolumn{1}{c|}{$\beta_{1,1}$}       & \multicolumn{1}{c|}{$\beta_{2,1}$}      & \multicolumn{1}{c|}{$\beta_{3,1}$}      & \multicolumn{1}{c|}{$\cdots$} & \multicolumn{1}{c|}{$\beta_{e-4,1}$}    & \multicolumn{1}{c|}{$\beta_{e-3,1}$}          & \multicolumn{1}{c|}{$\beta_{e-2,1}$}          & \multicolumn{1}{c|}{$\beta_{e-1,1}$}        & $\beta_{e,1}$       & Extremal case             & Reference                                                    \\ \hline\hline
				\multicolumn{1}{|c|}{\multirow{14}{*}{The first hierarchy}}  & \begin{tabular}[c]{@{}c@{}}$\propp[0,0]$\\ $(d\geq e+1)$\end{tabular}                 & \multicolumn{9}{c|}{$\leq p\binom{e+1}{p+1}$\quad\quad\quad\quad}                                                                                                                                                                                                                                                                                                                              & $\vmd[n]$                 & {\cite{MR3302629}}                                           \\ \cline{2-13} 
				\multicolumn{1}{|c|}{}                                      & \begin{tabular}[c]{@{}c@{}}$\propp[0,1]$\\ $(d\geq e+2)$\end{tabular}                 & \multicolumn{8}{c|}{$\leq p\binom{e+1}{p+1}-\binom{e}{p-1}$\quad}                                                                                                                                                                                                                                                                                         & \multirow{38}{*}{0} & $\dpv[n]$                 & {\cite{MR3302629}}                                           \\ \cline{2-10} \cline{12-13} 
				\multicolumn{1}{|c|}{}                                      & \begin{tabular}[c]{@{}c@{}}$\propp[0,2]$\\ $(d\geq e+3)$\end{tabular}                 & \multicolumn{8}{c|}{$\leq p\binom{e+1}{p+1}-\textcolor{blue}{2}\binom{e}{p-1}$\quad\quad\quad}                                                                                                                                                                                                                                                                                        &                     & $\acm[n,3]$               & {\cite{MR4645613}}                                           \\ \cline{2-10} \cline{12-13} 
				\multicolumn{1}{|c|}{}                                      & \begin{tabular}[c]{@{}c@{}}$\propp[0,3]$\\ $(d\geq e+4)$\end{tabular}                 & \multicolumn{7}{c|}{$\leq p\binom{e+1}{p+1}-\textcolor{blue}{3}\binom{e}{p-1}$\quad\quad}                                                                                                                                                                                                                                          & \multicolumn{1}{c|}{$\leq p\binom{e}{p+1}$} &                     & $\acm[n,4]$               & \multirow{35}{*}{\Cref{hierarchy_theorem_intro}}        \\ \cline{2-10} \cline{12-12}
				\multicolumn{1}{|c|}{}                                      & \begin{tabular}[c]{@{}c@{}}$\propp[0,4]$\\ $(d\geq e+5)$\end{tabular}                 & \multicolumn{6}{c|}{$\leq p\binom{e+1}{p+1}-\textcolor{blue}{4}\binom{e}{p-1}$}                                                                                                                                                                                          & \multicolumn{2}{c|}{$\leq p\binom{e}{p+1}$}                                                 &                     & $\acm[n,5]$               &                                                              \\ \cline{2-10} \cline{12-12}
				\multicolumn{1}{|c|}{}                                      & \begin{tabular}[c]{@{}c@{}}$\propp[0,5]$\\ $(d\geq e+6)$\end{tabular}                 & \multicolumn{5}{c|}{$\leq p\binom{e+1}{p+1}-\textcolor{blue}{5}\binom{e}{p-1}$}                                                                                                                                          & \multicolumn{3}{c|}{$\leq p\binom{e}{p+1}$}                                                                                                 &                     & $\acm[n,6]$               &                                                              \\ \cline{2-10} \cline{12-12}
				\multicolumn{1}{|c|}{}                                      & $\vdots$                                                                              & \multicolumn{8}{c|}{$\vdots$}                                                                                                                                                                                                                                                                                                                        &                     & $\vdots$                  &                                                              \\ \cline{2-10} \cline{12-12}
				\multicolumn{1}{|c|}{}                                      & \begin{tabular}[c]{@{}c@{}}$\propp[0,e-1]$\\ $(d\geq 2e)$\end{tabular}                & \multicolumn{2}{c|}{$\leq p\binom{e+1}{p+1}-\textcolor{blue}{(e-1)}\binom{e}{p-1}$}                  & \multicolumn{6}{c|}{$\leq p\binom{e}{p+1}$}                                                                                                                                                                                                                     &                     & $\acm[n,e]$               &                                                              \\ \cline{1-10} \cline{12-12} 
				\multicolumn{1}{|c|}{\multirow{10}{*}{The second hierarchy}} & \begin{tabular}[c]{@{}c@{}}$\propp[0,e]$\\ $(d\geq 2e+1)$\end{tabular}                & \multicolumn{8}{c|}{$\leq p\binom{e}{p+1}$\quad\quad\quad\quad\quad\quad\quad}                                                                                                                                                                                                                                                                                                          &                     & $\subseteq \vmd[n+1]$     &  \\ \cline{2-10} \cline{12-12} 
				\multicolumn{1}{|c|}{}                                      & \begin{tabular}[c]{@{}c@{}}$\propp[1,1]$\\ $(d\gg 0$, $\propa[1,1])$\end{tabular}     & \multicolumn{7}{c|}{$\leq p\binom{e}{p+1}-\binom{e-1}{p-1}$}                                                                                                                                                                                                                                           & \multicolumn{1}{c|}{\multirow{24}{*}{0}}    &                     & $\subseteq \dpv[n+1]$     &                  \\ \cline{2-9} \cline{12-12} 
				\multicolumn{1}{|c|}{}                                      & \begin{tabular}[c]{@{}c@{}}$\propp[1,2]$\\ $(d\gg 0$, $\propa[1,2])$\end{tabular}     & \multicolumn{7}{c|}{$\leq p\binom{e}{p+1}-\textcolor{blue}{2}\binom{e-1}{p-1}$\quad\quad\quad}                                                                                                                                                                                                                                          & \multicolumn{1}{c|}{}                       &                     & $\subseteq \acm[n+1,3]$   &           \\ \cline{2-9} \cline{12-12}
				\multicolumn{1}{|c|}{}                                      & \begin{tabular}[c]{@{}c@{}}$\propp[1,3]$\\ $(d\gg 0$, $\propa[1,3])$\end{tabular}     & \multicolumn{6}{c|}{$\leq p\binom{e}{p+1}-\textcolor{blue}{3}\binom{e-1}{p-1}$}                                                                                                                                                                                          & \multicolumn{1}{c|}{$\leq p\binom{e-1}{p+1}$} & \multicolumn{1}{c|}{}                       &                     & $\subseteq \acm[n+1,4]$   &                                                              \\ \cline{2-9} \cline{12-12}
				\multicolumn{1}{|c|}{}                                      & $\vdots$                                                                              & \multicolumn{7}{c|}{$\vdots$}                                                                                                                                                                                                                                                                          & \multicolumn{1}{c|}{}                       &                     & $\vdots$                  &                                                              \\ \cline{2-9} \cline{12-12}
				\multicolumn{1}{|c|}{}                                      & \begin{tabular}[c]{@{}c@{}}$\propp[1,e-2]$\\ $(d\gg 0$, $\propa[1,e-2])$\end{tabular} & \multicolumn{2}{c|}{$\leq p\binom{e}{p+1}-\textcolor{blue}{(e-2)}\binom{e-1}{p-1}$}                  & \multicolumn{5}{c|}{$\leq p\binom{e-1}{p+1}$}                                                                                                                                                                     & \multicolumn{1}{c|}{}                       &                     & $\subseteq \acm[n+1,e-1]$ &                                                              \\ \cline{1-9} \cline{12-12}
				\multicolumn{1}{|c|}{\multirow{7}{*}{The third hierarchy}}  & \begin{tabular}[c]{@{}c@{}}
					$\propp[1,e-1]$       \\
					$(d\gg 0$, $\propa[1,e-1])$
				\end{tabular} & \multicolumn{7}{c|}{$\leq p\binom{e-1}{p+1}$\quad\quad\quad\quad\quad\quad\quad}                                                                                                                                                                                                                                                          & \multicolumn{1}{c|}{}                       &                     & $\subseteq \vmd[n+2]$     &                                                              \\ \cline{2-9} \cline{12-12}
				\multicolumn{1}{|c|}{}                                      & $\propp[2,1]$                                                                         & \multicolumn{6}{c|}{$\leq p\binom{e-1}{p+1}-\binom{e-2}{p-1}$}                                                                                                                                                                                         & \multicolumn{1}{c|}{\multirow{13}{*}{0}}      & \multicolumn{1}{c|}{}                       &                     & $\subseteq \dpv[n+2]$     &                                                              \\ \cline{2-8} \cline{12-12}
				\multicolumn{1}{|c|}{}                                      & $\propp[2,2]$                                                                         & \multicolumn{6}{c|}{$\leq p\binom{e-1}{p+1}-\textcolor{blue}{2}\binom{e-2}{p-1}$\quad\quad\quad}                                                                                                                                                                                        & \multicolumn{1}{c|}{}                         & \multicolumn{1}{c|}{}                       &                     & $\subseteq \acm[n+2,3]$   &                                                              \\ \cline{2-8} \cline{12-12}
				\multicolumn{1}{|c|}{}                                      & $\propp[2,3]$                                                                         & \multicolumn{5}{c|}{$\leq p\binom{e-1}{p+1}-\textcolor{blue}{3}\binom{e-2}{p-1}$}                                                                                                                                        & \multicolumn{1}{c|}{$\leq p\binom{e-2}{p+1}$} & \multicolumn{1}{c|}{}                         & \multicolumn{1}{c|}{}                       &                     & $\subseteq \acm[n+2,4]$   &                                                              \\ \cline{2-8} \cline{12-12}
				\multicolumn{1}{|c|}{}                                      & $\vdots$                                                                              & \multicolumn{6}{c|}{$\vdots$}                                                                                                                                                                                                                          & \multicolumn{1}{c|}{}                         & \multicolumn{1}{c|}{}                       &                     & $\vdots$                  &                                                              \\ \cline{2-8} \cline{12-12}
				\multicolumn{1}{|c|}{}                                      & $\propp[2,e-3]$                                                                       & \multicolumn{2}{c|}{$\leq p\binom{e-1}{p+1}-\textcolor{blue}{(e-3)}\binom{e-2}{p-1}$}                & \multicolumn{4}{c|}{$\leq p\binom{e-2}{p+1}$}                                                                                                                     & \multicolumn{1}{c|}{}                         & \multicolumn{1}{c|}{}                       &                     & $\subseteq \acm[n+2,e-2]$ &                                                              \\ \cline{1-8} \cline{12-12}
				\multicolumn{1}{|c|}{\multirow{4}{*}{The 4-th hierarchy}}   & $\propp[2,e-2]$                                                                       & \multicolumn{6}{c|}{$\leq p\binom{e-2}{p+1}$}                                                                                                                                                                                                          & \multicolumn{1}{c|}{}                         & \multicolumn{1}{c|}{}                       &                     & $\subseteq \vmd[n+3]$     &                                                              \\ \cline{2-8} \cline{12-12}
				\multicolumn{1}{|c|}{}                                      & $\vdots$                                                                              & \multicolumn{5}{c|}{$\vdots$}                                                                                                                                                                          & \multicolumn{1}{c|}{\multirow{5}{*}{0}}       & \multicolumn{1}{c|}{}                         & \multicolumn{1}{c|}{}                       &                     & $\vdots$                  &                                                              \\ \cline{2-7} \cline{12-12}
				\multicolumn{1}{|c|}{}                                      & $\propp[3,e-4]$                                                                       & \multicolumn{2}{c|}{$\leq p\binom{e-2}{p+1}-\textcolor{blue}{(e-4)}\binom{e-3}{p-1}$}                & \multicolumn{3}{c|}{$\leq p\binom{e-3}{p+1}$}                                                                     & \multicolumn{1}{c|}{}                         & \multicolumn{1}{c|}{}                         & \multicolumn{1}{c|}{}                       &                     & $\subseteq \acm[n+3,e-3]$ &                                                              \\ \cline{1-7} \cline{12-12}
				\multicolumn{1}{|c|}{$\vdots$}                              & $\vdots$                                                                              & \multicolumn{5}{c|}{$\vdots$}                                                                                                                                                                          & \multicolumn{1}{c|}{}                         & \multicolumn{1}{c|}{}                         & \multicolumn{1}{c|}{}                       &                     & $\vdots$                  &                                                              \\ \cline{1-7} \cline{12-12}
				\multicolumn{1}{|c|}{The $e$-th hierarchy}                  & $\propp[e-2,2]$                                                                       & \multicolumn{1}{c|}{$\leq \binom{2}{2}$} & \multicolumn{1}{c|}{\multirow{2}{*}{0}} & \multicolumn{1}{c|}{\multirow{2}{*}{0}} & \multicolumn{1}{c|}{$\cdots$} & \multicolumn{1}{c|}{\multirow{2}{*}{0}} & \multicolumn{1}{c|}{}                         & \multicolumn{1}{c|}{}                         & \multicolumn{1}{c|}{}                       &                     & $\subseteq \vmd[n+e-1]$   &                                                              \\ \cline{1-3} \cline{6-6} \cline{12-12}
				\multicolumn{1}{|c|}{The last hierarchy}                    & $\propp[e-1,1]$                                                                       & \multicolumn{1}{c|}{0}                   & \multicolumn{1}{c|}{}                   & \multicolumn{1}{c|}{}                   & \multicolumn{1}{c|}{$\cdots$} & \multicolumn{1}{c|}{}                   & \multicolumn{1}{c|}{}                         & \multicolumn{1}{c|}{}                         & \multicolumn{1}{c|}{}                       &                     & all                       &                                                              \\ \hline
			\end{tabular}%
		}
		\vspace{3pt}
		\caption{The complete hierarchies}
		\label{table:complete hierarchies}
	\end{table}
	
	Now we give an application of this hierarchical structure.
	The $K_{p,1}$-theorem developed by M. Green \cite{MR739785}, Nagel--Pitteloud \cite{MR1291122}, and Brodmann--Schenzel \cite{MR2274517} is a fundamental theorem that states the properties of $\beta_{p,1}$ and their relation to the geometric properties.
	\begin{theorem}[{\cite{MR739785,MR1291122,MR2274517}}, The $K_{p,1}$-theorem]\label{Kp1_intro}
		Let $X\subseteq\p^r$ be a nondegenerate projective variety of dimension $n$, codimension $e\geq 2$, and degree $d$. Then the following hold:
		\begin{enumerate}
			\item[\textnormal{(i)}] $\beta_{p,1}(X)=0$ for $p\geq e+1$.
			\item[\textnormal{(ii)}] $\beta_{e,1}(X)\neq 0$ if and only if $\beta_{e,1}(X)=e$ and $X$ is a \vmd[n].
			\item[\textnormal{(iii)}] If $d\geq e+2$, then the following are equivalent:
			\begin{enumerate}
				\item[\textnormal{(1)}] $\beta_{e-1,1}(X)\neq 0$;
				\item[\textnormal{(2)}] $\beta_{e-1,1}(X)\in \{\binom{e+1}{2}-1,e-1\}$;
			\end{enumerate}
			\item[\textnormal{(iv)}] If $d\geq e+3$, then the following are equivalent:
			\begin{enumerate}
				\item[\textnormal{(1)}] $\beta_{e-1,1}(X)\neq 0$;
				\item[\textnormal{(2)}] $\beta_{e-1,1}(X)=e-1$;
				\item[\textnormal{(3)}] $X$ is a subvariety of a \vmd[n+1].
			\end{enumerate}
		\end{enumerate}
	\end{theorem}
	
	In particular, it shows there is a relation between the nonvanishing of $\beta_{e-1,1}(X)$ and the containment of $X$ in a \vmd[n+1], i.e., the nonvanishing of $\beta_{e-1,1}(X)$ has a geometric interpretation.
	Then naturally, one may raise the following question:
	\begin{question}
		What is the geometric meaning of $\beta_{p,1}(X)\neq 0$ for $p\leq e-2$?
	\end{question}
	The following generalized $K_{p,1}$-theorem demonstrates that the vanishing pattern of $\beta_{p,1}$ reflects the existence of a variety of minimal degree containing $X$ at each hierarchy. 
	Note that if one takes $k=-1$ in the following theorem, then it recovers the $K_{p,1}$-theorem (\Cref{Kp1_intro}).
	
	\begin{thmalphabetintro}[Generalized $K_{p,1}$-theorem]\label{Gen_Kp1_intro}
		Let $X\subseteq\p^r$ be a nondegenerate projective variety of dimension $n$, codimension $e\geq 2$, and degree $d$ satisfying $\propp[k,e-k]$ for some $-1\leq k\leq e-3$. Then the following hold.
		\begin{itemize}
			\item[\textnormal{(i)}] $\beta_{p,1}(X)=0$ for $p\geq e-k$;
			\item[\textnormal{(ii)}] $\beta_{e-k-1,1}(X)\neq 0$ if and only if $\beta_{e-k-1,1}(X)=e-k-1$ and $X$ is a subvariety of a \vmd[n+k+1];
			\item[\textnormal{(iii)}] If $X$ satisfies $\propp[k+1,1]$, then the following are equivalent:
			\begin{enumerate}
				\item[\textnormal{(1)}] $\beta_{e-k-2,1}(X)\neq 0$;
				\item[\textnormal{(2)}] $\beta_{e-k-2,1}(X)\in \{(e-k-2)\binom{e-k}{e-k-1}-\binom{e-k-1}{e-k-3},e-k-2\}$;
			\end{enumerate}
			\item[\textnormal{(iv)}] If $X$ satisfies $\propp[k+1,2]$, then the following are equivalent:
			\begin{enumerate}
				\item[\textnormal{(1)}] $\beta_{e-k-2,1}(X)\neq 0$;
				\item[\textnormal{(2)}] $\beta_{e-k-2,1}(X)=e-k-2$;
				\item[\textnormal{(3)}] $X$ is a subvariety of a \vmd[n+k+2].
			\end{enumerate}
		\end{itemize}
	\end{thmalphabetintro}
	
	In \Cref{hierarchy_upper_bound_intro}, the extremal varieties are ACM varieties contained in \vmd[n+1] for $3\leq m\leq e-1$.
	These varieties admit a precise characterization in terms of the ambient varieties containing them.
	
	\begin{thmalphabetintro}[Extremal varieties in the first hierarchy]\label{thm:ACM with degree=e+u_intro}
		Let $X \subseteq \p^r$ be a nondegenerate projective variety of dimension $n$, codimension $e$, and degree $d=e+1+m$ with $3 \leq m \leq e-1$. Then the following statements are equivalent:
		\begin{enumerate}
			\item[\textnormal{(i)}] $X$ is an ACM variety contained in \vmd[n+1].
			\item[\textnormal{(ii)}] One of the following holds:
			\begin{enumerate}
				\item[$(1)$] $e=4$, $d=8$, and $X$ lies on (a cone over) the Veronese surface as a divisor.
				\item[$(2)$] $X$ is contained in a unique $(n+1)$-fold rational normal scroll $Y = S(a_1 , a_2 , \ldots ,a_{n+1} )$ for some integers $0 \leq a_1 \leq a_2 \leq \ldots \leq a_{n+1}$ such that 
				\begin{enumerate}
					\item[$(\alpha)$] If $a_{n} \geq 1$ then the divisor class of $X$ is $2H+(1-e+m)F$  where $H$ and $F$ denote the hyperplane section and a ruling of $Y$, respectively.
					\item[$(\beta)$] If $a_{n} = 0$ then the divisor class of $X$ is $(e+1+m)F$  where $F \cong \p^n$ is the effective generator of the divisor class group of $Y$.
				\end{enumerate}
			\end{enumerate}
		\end{enumerate}
	\end{thmalphabetintro}
	
	The paper is organized as follows. Section 2 reviews the preliminary ingredients: the syzygies of finite sets of points, partial elimination ideals, and the behavior of graded Betti numbers under inner projections.
	We also collect several extremal criteria--maximal quadrics, the characterization of $\beta_{e-1,1}$ 
	for varieties contained in a VMD, and the gonality theorem--that will be used in later sections. Section 3 establishes the first hierarchy of upper bounds on linear syzygies and proves our main algebraic results, forming the foundation for the complete hierarchy developed later. Finally, Section 4 provides a geometric description of extremal varieties in the first hierarchy. We show that these extremal cases are realized either as divisors on rational normal scrolls or, in the codimension-four case, on a cone over the Veronese surface, and we explain the appearance of two eligible Betti tables when $m=3$.
	
	\section*{Acknowledgement}
	J. I. Han was partially supported by Basic Science Research Program through the National Research Foundation of Korea (NRF) funded by the Ministry of Education (2019R1A6A1A10073887) when he was at Korea Advanced Institute of Science and Technology (KAIST). J. I. Han is currently supported by a KIAS Individual Grant (MG101401) at Korea Institute for Advanced Study.
	S. Kwak was supported by the National Research Foundation of Korea (NRF) grant funded
	by the Korea government (MSIT) (No. RS-2024-00352592).
	W. Lee was supported by Basic Science Research Program through NRF funded by the Ministry of Education (no. RS-2023-00250807).
	The authors are grateful to Euisung Park for the valuable discussion. 
	
	\section{Preliminaries and known results}
	
	There may exist two ways to reduce the projective varieties: either by taking general hyperplane sections or by taking the images of projections. When we take general hyperplane sections, the results on syzygies of finite points are used importantly. The following is one of the most well-known results on finite points.
	
	\begin{theorem}[{\cite[Theorem 1]{MR959214}}]\label{thm_finite_pts_Np}
		Let $X\subseteq\p^r$ be a set of $(2r+1-p)$-points in the general position. Then $X$ satisfies the property $\N_{2,p}$.
	\end{theorem}
	
	Also, M. P. Cavaliere, M. E. Rossi, and G. Valla studied the syzygies of finite points in \cite{MR1110573}.
	The following proposition gives the relation between the Betti numbers in the diagonal position.
	
	\begin{proposition}[{\cite[Proposition 2.4]{MR1110573}}]
		Let $X\subseteq\p^r$ be a set of $d$ points in general position and $t$ be the initial degree of $I_X$. Then for $1\leq p\leq r$, the following are equivalent.
		\begin{enumerate}
			\item[\textnormal{(1)}] $X$ satisfies $\N_{2,p}$
			\item[\textnormal{(2)}] $\beta_{p,t}(X)=0$
			\item[\textnormal{(3)}] $\beta_{p+1,t-1}(X)\leq \binom{p+t-1}{p}\binom{r+t}{p+t}-d\binom{r}{p}$
		\end{enumerate}
	\end{proposition}
	
	Now we briefly review the partial elimination ideals and the effects of inner projections on the graded Betti numbers.
	Let $I\subseteq S=\Bbbk[x_0,\cdots,x_r]$ be a homogeneous ideal.
	Pick a point $q\in \p^r=\proj S$. By taking a coordinate change, we may assume $q=[1,0,\cdots,0]$.
	Denote $R=\Bbbk[x_1,\cdots,x_r]$.
	Then we define the $i$-th partial elimination ideal of $I$ with respect to $q$ as
	\[
	K_i(I,q)=\{f'\in R\mid \exists f\in I\text{ such that }f=f'x_0^i+(\text{lower $x_0$-degree terms})\}\cup\{0\}.
	\]
	One can easily see this is a homogeneous ideal of $R$.
	
	The effects of inner projections on the graded Betti numbers are studied in \cite{MR2946932,MR3302629}. The following result plays a crucial role to prove \Cref{upper_bound_second_floor_next_intro}.
	\begin{theorem}[{\cite[Theorem 3.1]{MR3302629}}]\label{inner_proj_Betti}
		Let $X\subseteq\p^{r}$ be a nondegenerate subscheme and pick a point $q\in X$. Denote the image of the inner projection of $X$ from $q$ as $X_q$. Set $t=\dim_\Bbbk(K_1(I_X,q))_1$ and $e=\codim_q(X,\p^r)$. Then
		\begin{align*}
			\beta_{p,1}(X) &\leq \beta_{p,1}(X_q)+\beta_{p-1,1}(X_q)+\binom{t}{p}\\
			\beta_{1,1}(X) &= \beta_{1,1}(X_q)+t.
		\end{align*}
	\end{theorem}
	
	The following proposition characterizes the projective varieties having the maximum number of quadrics.
	
	\begin{proposition}[{\cite[Corollary 4.4]{MR3334084}}]\label{quadric_ACM}
		Let $X\subseteq\p^r$ be a nondegenerate projective variety of codimension $e$ and degree $d=e+1+m$ where $0\leq m\leq e-1$. Then $\beta_{1,1}(X)=\binom{e+1}{2}-m$ if and only if it is ACM.
	\end{proposition}
	
	The number $\beta_{e-1,1}$ is known to have only two values unless $X$ is a variety of minimal degree or a del Pezzo variety.
	Specifically,
	\begin{theorem}[{\cite[Corollary 3.1]{MR4645613}}]\label{beta_e-1}
		If a nondegenerate projective variety $X\subseteq\p^r$ is contained in a variety of minimal degree as a divisor, then $\beta_{e-1,1}(X)=e-1$.
	\end{theorem}

	We also recall the gonality conjecture, which is now the gonality theorem, proved by L. Ein and R. Lazarsfeld. The optimal bound for the degree was given by J. Rathmann. It determines the length of the quadratic strand of sufficiently positive curves in terms of the gonality $\gon(C)$ of $C$.
	
	\begin{theorem}[{\cite[Gonality theorem]{MR3415069,rathmann2016effectiveboundgonalityconjecture}}]\label{gonality_conjecture}
		Let $C\subseteq\p^r$ be a linearly normal smooth curve of degree $d$ and genus $g\geq 2$. If $d\geq 4g-3$, then $\beta_{p,1}(C)\neq 0$ if and only if $1\leq p\leq r-\gon(C)$.
	\end{theorem}
	
	The following theorem classifies the cases violate the gonality theorem if we lower the degree from the optimal one.
	
	\begin{theorem}[{\cite[Theorem 1.1]{MR4808077}}]\label{recent_gonality}
		Let $C$ be a smooth projective curve of genus $g\geq 2$ and $\l$ be a line bundle on $C$ with $\deg \l\geq 4g-4$. Denote the Koszul cohomology of $C$ with respect to $\l$ as $K_{p,1}(C;\l)$. Then $$K_{p,1}(C;\l)\neq 0\text{ if and only if }1\leq p\leq \deg\l-g-\gon(C)$$ unless $\l=\omega_C^2$ and either $g=2$ or $C$ can be realized as a smooth quartic curve. In the exceptional case, $$K_{p,1}(C;\l)\neq 0\text{ if and only if }1\leq p\leq \deg\l-g-\gon(C)+1.$$
	\end{theorem}
	
	For general $\nu$-gonal curves, the condition on $\deg \mathcal{L}$ for the gonality conjecture can be improved significantly by the result of G. Farkas and M. Kemeny \cite{MR4008525}.
	It would be interesting to know how the gonality $\nu$ plays a role to the upper bounds of $\beta_{p,1}(X)$ in the nonvanishing interval.
	
	Finally, we point out that except few cases, the varieties are contained in at most one variety of minimal degree as a divisor.
	\begin{lemma}[{\cite[Lemma 4.6]{MR4645613}}]\label{lem:uniqueness of embedding scroll}
		Let $X \subseteq \p^r$ be a nondegenerate variety of dimension $n$ and codimension $e \geq 3$ and degree $d \geq e+3$. Then there exists at most one \vmd[n+1] containing $X$.
	\end{lemma}

	\section{Hierarchy on graded Betti numbers}
	First, we induce the sharp upper bound of Betti numbers in the quadratic strand.
	
	\begin{thmalphabetmaintext}[The first hierarchy]\label{thm_upper_bound_general_maintext}
		Let $X\subseteq\p^r$ be a nondegenerate projective variety of dimension $n$, codimension $e\geq 2$, and degree $d\geq e+1+m$ for an integer $0\leq m\leq e$. Then we have
		\begin{numcases}{\beta_{p,1}(X)\leq}
			p\binom{e+1}{p+1}-m\binom{e}{p-1} &for $1\leq p\leq e+1-m$\\
			p\binom{e}{p+1} &for $e+1-m\leq p$.
		\end{numcases}
		When $2\leq m\leq e-1$, the equality holds for all $p$ if and only if $X$ is an ACM divisor of degree $e+m+1$ in a variety of minimal degree, i.e., $X$ is an \acm[n,m+1].
	\end{thmalphabetmaintext}
	\begin{proof}
		Let $Y\subseteq\p^e$ be the general $e$-plane section of $X$. Then $Y$ is a set of $d$ points in general position in $\p^e$. Let $Z$ be a set of $e+m+1$ points contained in $Y$. Since $Z$ satisfies $\N_{2,e-m}$ by \Cref{thm_finite_pts_Np}, we have $\beta_{p,2}(Z)=0$ for all $p\leq e-m$.
		As it is well known that $\beta_{p,1}(Z)-\beta_{p-1,2}(Z)=p\binom{e+1}{p+1}-m\binom{e}{p-1}$ (cf. \cite[Proposition 1.6]{MR1110573}), we get $\beta_{p,1}(Z)=p\binom{e+1}{p+1}-m\binom{e}{p-1}$ for all $1\leq p\leq e+1-m$. Then we get inequality (3) as $\beta_{p,1}(Z)\geq \beta_{p,1}(Y)\geq \beta_{p,1}(X)$ by the Lefschetz theorem. For inequality (4), we assume $m\geq 3$ as the case $m=1$ follows from the $K_{p,1}$-theorem and the case $m=2$ follows from \Cref{beta_e-1}. Since $d\geq e+1+l$ for all $0\leq l\leq m$, by applying inequality (3), we get
		$\beta_{p,1}(X)\leq p\binom{e+1}{p+1}-l\binom{e}{p-1}$ for all $1\leq p\leq e+1-l$ and for all $0\leq l\leq m$. In particular, we have
		\begin{align*}
			\beta_{e+1-l,1}(X) &\leq (e+1-l)\binom{e+1}{e+2-l}-l\binom{e}{e-l}\\
			&= (e+1-l)\binom{e}{e+2-l}
		\end{align*}
		for all $0\leq l\leq m$. Therefore we get the inequality (4).
		
		Suppose $X$ is an extremal variety, i.e., $\beta_{p,1}(X)$ achieves the upper bounds for all $p\geq 1$. Then it is ACM of degree $e+1+m$ by \Cref{quadric_ACM}. It is contained in a variety of minimal degrees as a divisor due to the $K_{p,1}$-theorem.
		Conversely, let $X\subseteq\p^r$ be an ACM variety contained in a variety of minimal degree $Y$ as a divisor.
		It should have degree $d=e+1+m$ since if $d>e+1+m$, then $\beta_{1,1}(X)$ cannot achieve the upper bound.
		Then $\beta_{p,1}(X)=p\binom{e+1}{p+1}-m\binom{e}{p-1}$ for $1\leq p\leq e+1-m$ since the general point section $Z\subseteq\p^e$ of $X$ satisfies $\beta_{p,1}(Z)=p\binom{e+1}{p+1}-m\binom{e}{p-1}$ for $1\leq p\leq e+1-m$ as before.
		Also, from $\beta_{p,1}(X)\geq \beta_{p,1}(Y)$, we get $\beta_{p,1}(X)\geq p\binom{e}{p+1}$. In particular, $\beta_{p,1}(X)=p\binom{e}{p+1}$ for $e+1-m\leq p$ by the inequality (4).
		Therefore, $X$ corresponds to the extremal case.
	\end{proof}
	
	\begin{remark}
		Not all ACM varieties of degree $d=e+1+m$ achieve the upper bounds. Later in \Cref{upper_bound_genus}, one can see that there exists a nonhyperelliptic linearly normal smooth curve $C$ of genus $g\geq 2$ and degree $d\geq 4g-3$ that is not extremal. Also, when the codimension is greater than 2, there are two eligible Betti tables for ACM varieties of $d=e+4$ (cf. \Cref{betti_e+4}). Among those, only one achieves the upper bounds of \Cref{thm_upper_bound_general_maintext}.
	\end{remark}

	When $d\geq 2e+1$, we get the following upper bound as a corollary of \Cref{thm_upper_bound_general_maintext}.
	
	\begin{coralphabetmaintext}\label{upper_bound_second_floor}
		With the same notation as \Cref{thm_upper_bound_general_maintext}, when $d\geq 2e+1$, we have \[\beta_{p,1}(X)\leq p\binom{e}{p+1}\] for all $p\geq 1$.
		The extremal varieties are exactly the ones contained in a \vmd[n+1].
	\end{coralphabetmaintext}
	\begin{proof}
		The bound follows by taking $m=e$ in \Cref{thm_upper_bound_general_maintext}.
		If $\beta_{p,1}(X)=p\binom{e}{p+1}$ for all $p\geq 1$, then $\beta_{e-1,1}(X)=e-1$ so that $X$ is contained in a \vmd[n+1] by the $K_{p,1}$-theorem.
		Conversely, if $X$ is contained in a \vmd[n+1] $Y$, then $\beta_{p,1}(X)\geq \beta_{p,1}(Y)=p\binom{e}{p+1}$, so the equality holds.
	\end{proof}
	
	Note that this extends the classical theorem in Castelnuovo theory to higher linear syzygies in the sense that when $p=1$ we get $\beta_{1,1}(X)\leq \binom{e}{2}$. By the Castelnuovo theory, when $d\geq 2e+3$, the varieties having the maximal number of quadrics are always contained in a variety of minimal degree as a divisor. However, when $d=2e+1$ or $d=2e+2$, there are varieties not contained in a variety of minimal degree as a divisor while it has the maximal number of quadrics as \Cref{counterexample_2e+1} and \Cref{counterexample_2e+2} show.
	
	\begin{example}[$d=2e+1$, $\beta_{1,1}=\binom{e}{2}$, but not contained in a VMD]\label{counterexample_2e+1}
		Let $C\subseteq\p^6$ be the image of the map
		\[
		\p^1  \to      \p^6,\quad [s,t]  \mapsto [s^{11}, s^{10}t, s^9t^2, s^8t^3, s^7t^4, s^5t^6, t^{11}].
		\]
		Using Macaulay2 \cite{M2}, one can check $C$ has the following Betti table
		\[\begin{array}{c|cccccc}
			&0 &1 &2 &3 &4 &5 \\ \hline
			0&1 & & & & \\
			1& &10 &16 &9 & &\\
			2& &1 &14 &26 &20 &5
		\end{array}\]
		with $\deg C =11$. This curve is not contained in a variety of minimal degree due to the $K_{p,1}$-theorem.
	\end{example}
	
	\begin{example}[$d=2e+2$, $\beta_{1,1}=\binom{e}{2}$, but not contained in a VMD]\label{counterexample_2e+2}
		Let $C\subseteq\p^2=\proj \Bbbk[s,t,w]$ be a smooth plane curve defined by $s^6+t^6+w^6$ and $\nu_3:\p^2\to\p^9$ be the 3-uple Veronese embedding
		\begin{equation*} 
			[s:t:w] \mapsto [s^3 : s^2 t : s^2 w : st^2 : stw :sw^2 : t^3 : t^2 w : tw^2 : w^3].
		\end{equation*} Then $\nu_3(C)$ is a canonical curve with $\deg \nu_3(C)=18$ and $g(\nu_3(C))=10$.
		Then $h^0(\i_{\nu_3(C)}(2))=\binom{8}{2}$.
		Using Macaulay2, one can check $\nu_3(C)$ has the following Betti table.
		\[\begin{array}{c|ccccccccc}
			& 0 & 1  & 2   & 3   & 4   & 5   & 6   & 7  & 8 \\ \hline
			0 & 1 &    &     &     &     &     &     &    &  \\
			1 &   & 28 & 105 & 189 & 189 & 105 & 27  &    &   \\
			2 &   &    & 27  & 105 & 189 & 189 & 105 & 28 &   \\
			3 &   &    &     &     &     &     &     &    & 1
		\end{array}\]
		As $\beta_{7,1}(\nu_3(C))=0$, it is not contained in a variety of minimal degree as a divisor. Indeed, Petri's theorem shows that the ideal of a canonical curve of genus at least 4 is generated by quadrics unless it is trigonal or a smooth plane quintic. When the ideal is generated by quadrics, it is not contained in a variety of minimal degree as a divisor by the $K_{p,1}$-theorem.
	\end{example}
	
	For linearly normal smooth curves of high degree, we get the following sharp bounds as a corollary of \Cref{thm_upper_bound_general_maintext}.
	
	\begin{corollary}\label{upper_bound_genus}
		Let $C\subseteq\p\H^0(C, \mathcal L)$ be a linearly normal smooth curve of genus $g\ge 2$ with $\deg\mathcal L\ge 2g+1+\alpha$. Then
		\[
		\beta_{p,1}(C)\leq
		\begin{cases}
			p\binom{g+\alpha+1}{p+1}-g\binom{g+\alpha}{p-1} & \text{ for } 1\leq p\leq \alpha+1 \\
			
			p\binom{g+\alpha}{p+1} & \text { for } \alpha+1\le p \le g+\alpha-1.\\
		\end{cases}
		\]
		Assume furthermore $\alpha\geq 2g-4$. Then the equalities hold for all $p\ge 1$ if and only if $C$ is hyperelliptic of degree $2g+1+\alpha$. 
	\end{corollary}
	\begin{proof}
		By Riemann-Roch theorem, we have $h^0(C,\mathcal{L})=2g+1+\alpha-g+1=g+\alpha+2$, so $e=\codim C=g+\alpha$. Hence the inequality follows from \Cref{thm_upper_bound_general_maintext}. When $C$ is hyperelliptic, the Betti numbers are well known, and those achieve the upper bounds. Conversely, if the equalities hold for all $p\geq 1$, then the gonality must be two due to the gonality theorem (\Cref{gonality_conjecture}).
	\end{proof}
	
	The condition $\alpha\geq 2g-4$ in the second statement of the above corollary is optimal.
	Indeed, when $\alpha=2g-5$, there is a nonhyperelliptic curve whose Betti numbers achieve the upper bounds in \Cref{upper_bound_genus} as the following example shows.
	
	\begin{example}[Nonhyperelliptic extremal curve]
		Let $C$ be a smooth plane quartic in $\p^2$ and $\mathcal{L}=\omega_C^2$. The curve $C$ is embedded as a trigonal curve of degree $8$ and genus $3$ in $\p^5$.
		Then it is 3-regular and satisfies the property $\N_{2,1}$ thanks to Green's $(2g+1+p)$-theorem \cite[Theorem (4.a.1)]{MR739785}. Also, $\beta_{3,1}(C)\neq 0$ by \Cref{recent_gonality}.
		Then we get $\beta_{3,1}(C)=3$ by \Cref{beta_e-1}.
		This enables us to determine all Betti numbers as follows.
		\[\begin{array}{c|ccccc}
			&0 &1 &2 &3 &4 \\ \hline
			0&1 & & & \\
			1& &7 &8 &3 & \\
			2& & &6 &8 &3 
		\end{array}\]
		In particular, this curve achieves the upper bounds $\beta_{1,1}=7$, $\beta_{2,1}=8$, and $\beta_{3,1}=3$ in \Cref{thm_upper_bound_general_maintext} while it is not hyperelliptic.
	\end{example}
	
	We also get the sharp bounds to the canonical curves.
	\begin{corollary}
		Let $C\subseteq\p^{g-1}$ be a canonical curve. Then
		\[
		\beta_{p,1}(C)\leq p\binom{g-2}{p+1}
		\]
		for all $p$. When $g\geq 4$, the equalities hold for all $p$ if and only if $C$ is trigonal or a smooth plane quintic curve.
	\end{corollary}
	\begin{proof}
		The bound follows from \Cref{upper_bound_second_floor}. 
		The second statement follows from the self duality of the syzygies of canonical curves, Petri's theorem, and $K_{p,1}$-theorem.
		Indeed, the extremal cases occur if and only if $C$ is contained in a \vmd[2] by \Cref{upper_bound_second_floor}, i.e., $\beta_{g-3,1}(C)\neq 0$.
		This is equivalent to $\beta_{1,2}(C)\neq 0$, and this is the case if and only if $C$ is trigonal or a smooth plane quintic curve by Petri's theorem. 
	\end{proof}
	
	Note that the following is true:
	\begin{itemize}
		\item When $d\geq e+1$, we have $\beta_{p,1}(X)\leq p\binom{e+1}{p+1}$ and the extremal cases are the varieties of minimal degree.
		\item When $d\geq e+1$, we have $\beta_{p,1}(X)\leq p\binom{e+1}{p+1}-\binom{e}{p-1}$ unless it is a variety of minimal degree, and the extremal cases are del Pezzo varieties.
		\item When $d\geq 2e+1$, we have $\beta_{p,1}(X)\leq p\binom{e}{p+1}$ and the extremal cases are the varieties contained in varieties of minimal degree as a divisor.
	\end{itemize}
	
	Suppose $d\geq 2e+3$ and $X$ is not a divisor of a variety of minimal degree. Then we have $\beta_{1,1}(X)\leq \binom{e}{2}-1$ from the Castelnuovo theory. It is natural to ask whether there exists an extension of this bound to higher syzygies.
	Indeed, we can extend this bound to syzygies when $d
	\gg 0$.
	
	\begin{thmalphabetmaintext}\label{next_upper_bound}\label{second_row_in_second_hierarchy}
		Let $X\subseteq\p^r$	 be a nondegenerate projective variety of dimension $n$, codimension $e\geq 2$, and degree $d$ that is not contained in a \vmd[n+1]. If $d\gg 0$, then 
		\[\beta_{p,1}(X)\leq p\binom{e}{p+1}-\binom{e-1}{p-1}\]
		for all $p\geq 1$. Furthermore, the extremal varieties are the ones contained in a \dpv[n+1] when $e\geq 3$.
	\end{thmalphabetmaintext}
	\begin{proof}
		When $e\leq 3$, it is trivial. Assume $e\geq 4$.
		As $X$ is not contained in a \vmd[n+1], we may assume $\beta_{1,1}(X)\leq \binom{e}{2}-1$. If $\beta_{1,1}(X)=\binom{e}{2}-1$, then $X$ is contained in a del Pezzo variety of dimension $n+1$ so that $\beta_{p,1}(X)=p\binom{e}{p+1}-\binom{e-1}{p-1}$ for all $p\geq 1$. Hence we assume $\beta_{1,1}(X)\leq \binom{e}{2}-2$.
		Suppose $d>2^{\binom{e}{2}-2}$. We use the induction on $e$.
		Let $W$ be the scheme defined by the ideal generated by $(I_X)_2$ and $Z_1,\cdots,Z_k$ be the irreducible components of $W_{\text{red}}$. Then $\sum_i\deg(Z_i)\leq 2^{\binom{e}{2}-2}$, so $X$ is not an irreducible component of $W_{\text{red}}$.
		This implies that for a general point $q$ of $X$, we have $\dim T_q W > \dim X$.
		By \Cref{inner_proj_Betti}, for the general point $q\in X$, we get $\beta_{p,1}(X)\leq \beta_{p,1}(X_q)+\beta_{p-1,1}(X_q)+\binom{t}{p}$ where $t=\dim_\Bbbk (K_1(I_X,q))_1$. As the elements of $(K_1(I_X,q))_1$ are the linear forms defining $T_qW$, we get $\dim_\Bbbk (K_1(I_X,q))_1=\codim T_qW<\codim X=e$. Hence $\beta_{p,1}(X)\leq \beta_{p,1}(X_q)+\beta_{p-1,1}(X_q)+\binom{e-1}{p}$. Note that $d-1>2^{\binom{e-1}{2}-2}$. Also, $X_q\subseteq\p^{r-1}$ is not contained a \vmd[n+1] by \cite[Theorem 4.1]{QPPY}. Hence $\beta_{p,1}(X_q)\leq p\binom{e-1}{p+1}-\binom{e-2}{p-1}$ for all $p\geq 1$.	Therefore we get $\beta_{p,1}(X)\leq p\binom{e-1}{p+1}-\binom{e-2}{p-1}+(p-1)\binom{e-1}{p}-\binom{e-2}{p-2}+\binom{e-1}{p}=p\binom{e}{p+1}-\binom{e-1}{p-1}$.
	\end{proof}
	
	The following gives an example of an extremal variety satisfying the condition of \Cref{second_row_in_second_hierarchy}.
	
	\begin{example}
		Let $C\subseteq\p^2=\proj \Bbbk[s,t,w]$ be a smooth plane octic curve defined by 
		$$s^8+t^8-w^8=0,$$
		and let $\nu_3:\p^2\to\p^9$ be the 3-uple Veronese embedding. Then $\nu_3(C)$ has degree $24$ and genus $21$, and $\nu_3(\mathbb{P}^2)$ is a del Pezzo surface of degree $9$. Consider the four points
		$$P_1=[1:0:1],\ P_2=[1:0:-1],\ P_3=[0:1:1],\  \mbox{and}\  P_4=[0:1:-1]$$
		on $C$. Their images under $\nu_3$ are 
		\begin{equation*} 
			\begin{cases}
				\nu_3(P_1)=[1:0:1:0:0:1:0:0:0:1],\\
				\nu_3(P_2)=[1:0:-1:0:0:1:0:0:0:-1],\\
				\nu_3(P_3)=[0:0:0:0:0:0:1:1:1:1], \ \mbox{and}\\
				\nu_3(P_4)=[0:0:0:0:0:0:1:-1:1:-1].
			\end{cases}
		\end{equation*}
		Let $\Lambda \cong \mathbb{P}^3$ be the linear span of these four linearly independent points, and consider the linear projection $\pi : \mathbb{P}^9 \rightarrow \mathbb{P}^5$ from  $\Lambda$. Let $f:= \pi|_{\nu_3(C)}$ and $g:= \pi|_{\nu_3(\mathbb{P}^2)}$ be the retrictions to $\nu_3(C)$ and $\nu_3(\mathbb{P}^2)$, respectively. Set $D:=f(\nu_3(C))$ and $S:=g(\nu_3(\mathbb{P}^2))$. Then $D \subset \mathbb{P}^5$ is the curve of degree $20$ and $S \subset \mathbb{P}^5$ is the del Pezzo surface of degree $5$. Moreover, $D$ has the Betti table
		\[\begin{array}{c|cccccc}
			& 0 & 1  & 2   & 3   & 4   & 5    \\ \hline
			0 & 1 &    &     &     &     &   \\
			1 &   & 5  & 5   &     &     &      \\
			2 &   &    &     & 1   &     &     \\
			3 &   &    &     &     &     &   \\
			4 &   &    &     &     &     &   \\
			5 &   & 6  & 15  & 12  & 3   &     \\
			6 &   &    & 7   & 18  & 15  & 4    \\
		\end{array}\]
		by Macaulay2.
		Note that $D$ is not contained in a surface of minimal degree as a divisor by the $K_{p,1}$-theorem. Indeed, $D$ is a extremal case described in Theorem \ref{next_upper_bound}. That is, $D$ lies on the del Pezzo surface $S$ cut out by the five quadrics in the table.
	\end{example}
	
	\Cref{next_upper_bound} suggests the existence of further hierarchies. Indeed, the following theorem establishes hierarchies in full generality.
	
	\begin{thmalphabetmaintext}[The complete hierarchies]\label{hierarchy_theorem_maintext}
		Let $X\subseteq\p^r$ be a nondegenerate projective variety of dimension $n$, codimension $e\geq 2$, and degree $d$. For some $0\leq k\leq e-1$ and $0\leq m\leq e-k$, if X satisfies $\propp[k,m]$, then
		\begin{numcases}{\beta_{p,1}(X)\leq}
			p\binom{e+1-k}{p+1}-m\binom{e-k}{p-1} &for $1\leq p\leq e+1-m-k$ \nonumber\\
			p\binom{e-k}{p+1} &for $e+1-m-k\leq p$. \nonumber
		\end{numcases}
		Furthermore, the extremal cases are determined as follows except $(k,m)=(e-1,1)$.
		\begin{itemize}
			\item[\textnormal{(i)}] When $m=0$, the equalities hold for all $p\geq 1$ if and only if $X$ is contained in a \vmd[n+k].
			\item[\textnormal{(ii)}] When $m=1$, the equalities hold for all $p\geq 1$ if and only if $X$ is contained in a \dpv[n+k].
			\item[\textnormal{(iii)}] When $2\leq m\leq e-k-1$, the equalities hold for all $p\geq 1$ if and only if $X$ is contained in an \acm[n+k,m+1].
			\item[\textnormal{(iv)}] When $m=e-k$, the equalities hold for all $p\geq 1$ if and only if $X$ is contained in a \vmd[n+k+1].
		\end{itemize}
		When $(k,m)=(e-1,1)$, all varieties (satisfying \propp[e-1,1]) are vacuously extremal.
	\end{thmalphabetmaintext}
	\begin{proof}
		When $k=0$, it follows from \Cref{thm_upper_bound_general_maintext}.	
		If $k=e-1$, then $m=0$ or $1$. When $m=1$, then the bounds are trivial, and when $m=0$, the only non-vacuous bound is $\beta_{1,1}(X)\leq 1$.
		For this case, as $X$ satisfies \propp[e-1,0], it satisfies \propa[e-2,d_{e-2}], i.e., \propa[e-2,2]. 
		Suppose $\beta_{1,1}(X)\geq 2$. Then there are two distinct nondegenerate quadrics $Q_1,Q_2\in (I_X)_2$. 
		Then irreducible components of the scheme $B$ defined by $Q_1,Q_2$ have degree at most 4. As $B$ is equidimensional of codimension two, this contradicts to that $X$ satisfies \propa[e-2,2].
		Hence the bounds hold when $k=e-1$.
		To see the extremal varieties, when $(k,m)=(e-1,0)$, as the only non-vacuous bound is $\beta_{1,1}(X)=1$, the extremal varieties are the ones contained in a \vmd[n+e-1].
		When $(k,m)=(e-1,1)$, as the bounds are vacuous, the extremal varieties are all varieties (satisfying \propp[e-1,1]).
		
		From now on, we assume $1\leq k\leq e-2$. 
		For these cases, we may assume $\beta_{1,1}(X)\geq 2$, since if not, we have $\beta_{1,1}(X)\leq 1$ and $\beta_{p,1}(X)=0$ for $p\geq 2$ which satisfy the bounds.
		Let $W$ be the scheme defined by $(I_X)_2$ and $Y$ be the irreducible component of $W$ containing $X$. Then $(I_X)_2\supseteq (I_Y)_2\supseteq (I_W)_2$. As $(I_X)_2=(I_W)_2$, we have $(I_X)_2=(I_Y)_2$. In particular, $\beta_{p,1}(X)=\beta_{p,1}(Y)$ for all $p\geq 1$.
		By the given condition, $X$ satisfies $\propa[0,d_0],\propa[1,d_1],\cdots,\propa[k-1,d_{k-1}]$. 
		
		We claim $e_Y:=\codim Y \leq e-k$. Suppose, on the contrary, that $e_Y\geq e-k+1$. Since $X$ satisfies the property $\propa[e-e_Y,d_{e-e_Y}]$, we obtain $d_Y:=\deg Y>e_Y+d_{e-e_Y}=2^{e_Y}$. Write the irreducible decomposition of $W$ as
		$$W = Y \cup W_1 \cup \cdots \cup W_{t_1}.$$  
		Choose two distinct quadratic forms $Q_1,Q_2\in (I_Y)_2\backslash \{0\}$ and denote $Z_2$ the scheme defined by $Q_1,Q_2$.
		Let 
		$$Z_2 = Z_{2,1} \cup Z_{2,2} \cup \cdots \cup Z_{2,t_2}$$ 
		be the irreducible decomposition of $Z_2$.
		Note that all of $Z_{2,1},Z_{2,2},\cdots,Z_{2,t_2}$ have codimension two.
		As $Y$ is irreducible and contained in $Z_2$, we may assume $Y\subseteq Z_{2,1}$. Thus $\deg Z_{2,1} \leq 2^2$.
		If $e_Y=2$, then $Y=Z_{2,1}$ so that $d_Y\leq 2^2$. This contradicts to $d_Y>2^{e_Y}$.
		If $e_Y>2$, then there is a nonzero quadratic form in $(I_Y)_2$ that does not vanish identically on $Z_{2,1}$. Otherwise, $Z_{2,1}$ would be an irreducible component of $W$ containing $Y$ properly, contradicting that $Y$ is the component lying in the base scheme of $(I_Y)_2$. Denote the above quadratic form as $Q_3$.
		Let $Z_3$ be the scheme defined by $I_{Z_{2,1}}+\langle Q_3\rangle$, and $Z_3 = Z_{3,1} \cup Z_{3,2}\cup\cdots \cup Z_{3,t_3}$ be the irreducible decomposition of $Z_3$. Similarly, we may say  $Y\subseteq Z_{3,1}$ and $\deg Z_{3,1} \leq 2^3$. Then repeat the above argument until we get $Z_{e_Y}$ which is defined by $I_{Z_{e_Y-1,1}}+\langle Q_{e_Y}\rangle$, and the irreducible decomposition
		$$Z_{e_Y} = Z_{e_Y,1}\cup Z_{e_Y,2}\cup  \cdots \cup Z_{e_Y,t_{e_Y}}$$
		of $Z_{e_Y}$. As before, it holds that $\deg Z_{e_Y,1} \leq 2^{e_Y}$ and we may assume $Y\subseteq Z_{e_Y,1}$. On the other hand, since $Z_{e_Y,1}$ has codimension $e_Y$ and is irreducible, we have the equality $Y = Z_{e_Y,1}$. This contradicts to the assumption $d_Y>2^{e_Y}$. Therefore the claim is showed.
		
		Then we divide into two cases.
		First, we assume $e_Y=e-k$.
		As $X$ satisfies $\propa[k,m]$, we have $d_Y\geq e_Y+1+m$.
		Then by applying \Cref{thm_upper_bound_general_maintext} to $Y$, we get
		\begin{numcases}{\beta_{p,1}(Y)\leq}
			p\binom{e_Y+1}{p+1}-m\binom{e_Y}{p-1} &for $1\leq p\leq e_Y+1-m$ \nonumber \\
			p\binom{e_Y}{p+1} &for $e_Y+1-m\leq p$. \nonumber
		\end{numcases}
		Hence the case $e_Y=e-k$ is proved.
		
		Next, we consider the remaining case where $e_Y\leq e-k-1$. In this case, we have $\beta_{p,1}(Y)\leq p\binom{e_Y+1}{p+1}\leq p\binom{e-k}{p+1}$ by applying \Cref{thm_upper_bound_general_maintext} to $Y$ with $m=0$. 
		As $p\binom{e-k}{p+1}\leq p\binom{e+1-k}{p+1}-m\binom{e-k}{p-1}$ for any $1\leq p\leq e+1-m-k$ and any $0\leq m\leq e-k$, we get $\beta_{p,1}(X)\leq p\binom{e+1-k}{p+1}-m\binom{e-k}{p-1}$ for $1\leq p\leq e+1-m-k$. Note that the part $\beta_{p,1}(X)\leq p\binom{e-k}{p+1}$ for $e+1-m-k\leq p$ is already obtained.
		
		For the second statement, note that we have already determined the extremal varieties when $k=0$ or $k=e-1$.
		Hence we consider the case $1\leq k\leq e-2$.
		Suppose $X$ is contained in \vmd[n+k] ($m=0$), \dpv[n+k] ($m=1$), \acm[n+k,m+1] ($2\leq m\leq e-k-1$), or \vmd[n+k+1] ($m=e-k$), and denote it by $Z$.
		Then $\beta_{p,1}(X)\geq \beta_{p,1}(Z)$.
		Note that $\beta_{p,1}(Z)=p\binom{e+1-k}{p+1}-m\binom{e-k}{p-1}$ for $1\leq p\leq e+1-m-k$ and $\beta_{p,1}(Z)=p\binom{e-k}{p+1}$ for $e+1-m-k\leq p$.
		Hence the upper bounds of $\beta_{p,1}$ that we have just showed implies the equalities.
		Conversely, suppose the equality of the bounds hold for all $p\geq 1$. When $k=0$, the extremal cases are already determined above.
		When $k\geq 1$ and $e_Y=e-k$,  the variety $Y$ should be
		\begin{itemize}
			\item a \vmd[n+k] when $m=0$;
			\item a \dpv[n+k] when $m=1$;
			\item an \acm[n+k,m+1] when $2\leq m\leq e-k-1$;
			\item contained in \vmd[n+k+1] when $m=e-k$
		\end{itemize}
		by the extremal cases of \Cref{thm_upper_bound_general_maintext} and \Cref{upper_bound_second_floor}.
		Now assume $k\geq 1$ and $e_Y\leq e-k-1$. Then $\beta_{1,1}(Y)\leq \binom{e_Y+1}{2}\leq \binom{e-k}{2}$.
		As $\beta_{1,1}(X)=\binom{e+1-k}{2}-m$, we get $\binom{e+1-k}{2}-m\leq \binom{e-k}{2}$.
		Thus $m\geq e-k$, so $m=e-k$.
		Therefore $\beta_{e-k-1,1}(X)=(e-k-1)\binom{e-k}{e-k-1+1}=e-k-1\neq 0$.
		Note that if $e_Y\leq e-k-2$, then $\beta_{e-k-1,1}(Y)=0$ by the $K_{p,1}$-theorem which contradicts to $\beta_{e-k-1,1}(X)\neq 0$. Thus $e_Y=e-k-1$.
		As $\beta_{e-k-1,1}(Y)\neq 0$, by the $K_{p,1}$-theorem again, $Y$ is a \vmd[n+k+1] and the result follows.
	\end{proof}
	
	Using \Cref{hierarchy_theorem_maintext}, we prove the generalized $K_{p,1}$-theorem.
	
	\begin{thmalphabetmaintext}
		[Generalized $K_{p,1}$-theorem]\label{gen_Kp1}
		Let $X\subseteq\p^r$ be a nondegenerate projective variety of dimension $n$, codimension $e\geq 2$, and degree $d$ satisfying $\propp[k,e-k]$ for some $-1\leq k\leq e-3$. Then the following hold.
		\begin{itemize}
			\item[\textnormal{(i)}] $\beta_{p,1}(X)=0$ for $p\geq e-k$;
			\item[\textnormal{(ii)}] If $\beta_{e-k-1,1}(X)\neq 0$, then $\beta_{e-k-1,1}(X)=e-k-1$ and $X$ is a subvariety of a \vmd[n+k+1];
			\item[\textnormal{(iii)}] If $X$ satisfies $\propp[k+1,1]$, then the following are equivalent:
			\begin{enumerate}
				\item[\textnormal{(1)}] $\beta_{e-k-2,1}(X)\neq 0$;
				\item[\textnormal{(2)}] $\beta_{e-k-2,1}(X)\in \{(e-k-2)\binom{e-k}{e-k-1}-\binom{e-k-1}{e-k-3},e-k-2\}$;
			\end{enumerate}
			\item[\textnormal{(iv)}] If $X$ satisfies $\propp[k+1,2]$, then the following are equivalent:
			\begin{enumerate}
				\item[\textnormal{(1)}] $\beta_{e-k-2,1}(X)\neq 0$;
				\item[\textnormal{(2)}] $\beta_{e-k-2,1}(X)=e-k-2$;
				\item[\textnormal{(3)}] $X$ is a subvariety of a \vmd[n+k+2].
			\end{enumerate}
		\end{itemize}
	\end{thmalphabetmaintext}
	\begin{proof}
		Let $W$ be the scheme defined by $(I_X)_2$ and $Y$ be the irreducible component of $W$ containing $X$. As in the proof of \Cref{hierarchy_theorem_maintext}, the quadratic strands of $X$ is the same as that of $Y$ and we denote $e_Y:=\codim Y$.
		The statement (i) directly follows from \Cref{hierarchy_theorem_maintext}.
		
		To show (ii), assume $\beta_{e-k-1,1}(X)\neq 0$. Then $\beta_{e-k-1,1}(Y)\neq 0$. By applying the usual $K_{p,1}$-theorem to $Y$, we get $e_Y\geq e-k-1$.
		Recall that in the proof of \Cref{hierarchy_theorem_maintext}, we proved that if $X$ satisfies $\propa[0,d_0],\propa[1,d_1],\cdots,\propa[k-1,d_{k-1}]$, then $e_Y\leq e-k$.
		Hence $e_Y=e-k-1$ or $e-k$.
		When $e_Y=e-k-1$, the usual $K_{p,1}$-theorem implies $Y$ is a \vmd[n+k+1].
		When $e_Y=e-k$, as $X$ satisfies \propa[k,e-k], we have $\deg Y\geq e_Y+e-k+1\geq e_Y+4$, so the usual $K_{p,1}$-theorem shows $Y$ is contained in a \vmd[n+k+1].
		In any case, $X$ is contained in a \vmd[n+k+1] so that $\beta_{e-k-1,1}(X)\geq e-k-1$.
		By applying $(k,m)=(k,e-k)$ to \Cref{hierarchy_theorem_maintext}, we have $\beta_{e-k-1,1}(X)\leq e-k-1$. Hence we get the statement (ii).
		
		To show (iii), assume $X$ satisfies $\propp[k+1,1]$. As $X$ satisfies $\propa[0,d_0],\propa[1,d_1],\cdots,\propa[k,d_{k}]$, we have $e_Y\leq e-k-1$ as before.
		Also, as $X$ satisfies $\propa[k+1,1]$, the variety $Y$ is not a \vmd[n+k+1].
		Suppose $\beta_{e-k-2,1}(Y)\neq 0$.
		If $e_Y=e-k-1$, then by applying \Cref{Kp1_intro} to $Y$, we get $\beta_{e-k-2,1}(Y)\in \{(e-k-2)\binom{e-k}{e-k-1}-\binom{e-k-1}{e-k-3},e-k-2\}$.
		If $e_Y=e-k-2$, then $Y$ is a \vmd[n+k+2] and $\beta_{e-k-2,1}(Y)=e-k-2$.
		If $e_Y\leq e-k-3$, then we always have $\beta_{e-k-2,1}(Y)=0$, a contradiction.
		
		To show (iv), assume $X$ satisfies $\propp[k+1,2]$. As before, we have $e_Y\leq e-k-1$ and we divide into two cases.
		First, if $e_Y=e-k-1$, then by the condition $\propp[k+1,2]$, we have $\deg Y\geq e_Y+3$. Thus the result follows by applying \Cref{Kp1_intro} to $Y$.
		If $e_Y\leq e-k-2$, then whenever $\beta_{e-k-2,1}(Y)\neq 0$, $Y$ should be a \vmd[n+k+2]. Hence we get (iv.1)$\Rightarrow$(iv.3). 
		Also by \Cref{hierarchy_theorem_maintext}, we have $\beta_{e-k-2,1}(X)\leq e-k-2$. As (iv.3) implies $\beta_{e-k-2,1}(X)\geq e-k-2$, we get (iv.3)$\Rightarrow$(iv.2). Finally, (iv.2)$\Rightarrow$(iv.1) is trivial.
	\end{proof}

	\section{Extremal varieties in the first hierarchy}
	Let $X\subseteq\p^r$ be a nondegenerate projective variety of codimension $e$ and degree $d \geq e+1+m$ for $0\leq m\leq e-1$. 	
	This section is devoted to a geometric description of ACM varieties of degree $d=e+1+m$ with $3 \leq m \leq e-1$ contained in a variety of minimal degree to describe extremal varieties in the first hierarchy. 
	We only consider when $3 \leq m \leq e-1$ since the case $m=2$ is described in \cite{MR4645613}.
	If the Betti numbers $\beta_{p,1}(X)$ of $X$ attain the upper bounds in \Cref{thm_upper_bound_general_maintext} for all $p$, then $X$ must be an ACM variety of degree $e+1+m$. When $m=0,1,$ and $2$, the converse also holds. In contrast, when $m=3$, there exist ACM varieties of degree $e+4$ that fail to be extremal since in this case two distinct eligible Betti tables occur. 
	Recall that a set of seven points in general position in $\p^3$ admits exactly two eligible Betti tables (cf. \cite[Theorem 2.8]{MR2103875}). In general, we have 
	\begin{proposition}\label{betti_e+4}
		Let $X\subseteq\p^r$ be a nondegenerate projective variety of codimension $e \geq 3$ and degree $d=e+4$. Then $X$ is ACM if and only if its minimal free resolution has one of the following two Betti tables:
		\vspace{0.3cm}
		{\Small
			\[\begin{array}{c|cccccc}
				&0 &1 &\cdots &e-2 &e-1 &e \\ \hline
				0&1 & & & & & \\
				1& &\beta_{1,1} &\cdots &\beta_{e-2,1} & & \\
				2& & & &\binom{e-1}{2} &2e &3 
			\end{array}
			\mspace{60mu}
			\begin{array}{c|cccccc}
				&0 &1 &\cdots &e-2 &e-1 &e \\ \hline
				0&1 & & & & & \\
				1& &\beta_{1,1} &\cdots &\beta_{e-2,1} &e-1 & \\
				2& & & &\binom{e}{2} &2e &3 
			\end{array}\]
		}
		\vspace{0.3cm}
		\noindent where $\beta_{p,1}=p\binom{e+1}{p+1}-3\binom{e}{p-1}$ for all $p=1,\cdots,e-2$ in both tables.
	\end{proposition}
	\begin{proof}
		As in the proof of Theorem \ref{thm_upper_bound_general_maintext}, by applying $d=e+4$ to a general zero-dimensional section $Y$ of $X$, we obtain $\beta_{p,1}(X)-\beta_{p-1,2}(X)=p\binom{e+1}{p+1}-3\binom{e}{p-1}$ for all $p \geq 1$. Therefore the remaining task is to determine the values of $\beta_{p,1}(X)$, which are explicitly governed by \Cref{thm_upper_bound_general_maintext} and \Cref{beta_e-1}.
	\end{proof}
	
	Observe that among the two Betti tables above, only the second one attains the upper bounds for all $p$. Geometrically, this dichotomy reflects whether the variety lies as a divisor on a variety of minimal degree, as explained by Green's $K_{p,1}$-theorem. Hence ACM varieties of degree $e+4$ occur at the smallest degree for which such varieties may be either extremal or non-extremal. \\
	
	Let $X\subseteq\p^r$ be an ACM variety of degree $d=e+1+m$ for $3 \leq m \leq e-1$ that is contained in a \vmd[n+1].
	As a first step, we analyze the case when the variety is a curve. 
	
	\begin{proposition}\label{prop:unique embedding surface of minimal degree}
		Let $C \subset \p^{e+1}$ be a linearly normal curve of codimension $e$ and degree $d=e+1+m$ with $3 \leq m \leq e-1$. Suppose that $\beta_{e-1,1} (C) = e-1$. Then the arithmetic genus of $C$ is $m$, and $C$ is one of the following:
		\begin{enumerate}
			\item[\rm (1)] $e=4$, $d=8$, and $C$ is contained in the Veronese surface.
			\item[\rm (2)] $C$ is hyperelliptic and contained in the rational normal surface scroll $S(a,e-a)$ for some $\frac{e-1-m}{2} \leq a \leq \frac{e}{2}$. Furthermore, $C \equiv 2H+(1-e+m)F$ if $a > 0$, and $C \equiv 2eF$ if $a=0$, where $H$ and $F$ denote the hyperplane section and a ruling line of $S(a,e-a)$, respectively.
			\item[\rm (3)]  $C$ is nonhyperelliptic and contained in the singular rational normal surface scroll $S(0,e)$. Furthermore, $C \equiv (e+1+m)F$ where $F \cong \p^1$ is the effective generator of the divisor class group of $S(0,e)$.
		\end{enumerate}
		Moreover, the above rational normal surface scroll is the unique surface of minimal degree containing $C$.
	\end{proposition}
	
	\begin{proof}
		By the $K_{p,1}$-theorem and \Cref{lem:uniqueness of embedding scroll}, the curve $C$ is contained in a unique surface $S$ of minimal degree which is the Veronese surface or a rational normal surface scroll.
		By Riemann-Roch theorem, the arithmetic genus of $C$ is $m$. We divide into three cases.
		
		(1) Suppose $S$ is the Veronese surface in $\p^5$. 
		As $e=4$, we have $m=3$, so $C$ is the Veronese reembedding of a plane quartic curve $D\subseteq\p^2$.
		As $C$ has the maximal arithmetic genus, it is ACM by \cite{zbMATH02692308}.
		
		(2) Suppose $S$ is a rational normal surface scroll $S(a,e-a)$ ($0\leq a\leq \frac{e}{2}$) and $C$ is hyperelliptic. Then there exists a unique $f:C\to\p^1$ of degree 2. Let $A=f^*\o_{\p^1}(1)$ and write $\o_C(1)=A^l\otimes B$ where $B$ is an effective divisor of degree $b$ on $C$ such that $h^0(C,B)>0$ and $h^0(C,A^{-1}\otimes B)=0$. By considering degrees, we get $e+1+m=2l+b$. By \cite[Theorem 3.1.(2)]{MR2559684}, either $b=0,1$ and $l+b\geq m+1$, or $2\leq b\leq m+1$ and $l+b\geq m+2$ hold.
		By \cite[Corollary 3.3]{MR2559684}, we have $C\subseteq S(l+b-1-m,l)$.
		Hence $a=l+b-1-m$ and $a\geq \frac{e-1-m}{2}$.
		Write $C\equiv \alpha H+\beta F$. Here, when $a=0$, we regard $\alpha=0$ and $H:=0$.
		Then $\alpha\geq 0$.
		If $a\geq 1$, then $\alpha=2$ by the construction of $S(l+b-1-m,l)$.
		As $e+1+m=\deg C=(2H+\beta F)\cdot H=2e+\beta$, we get $\beta=1-e+m$.
		If $a=0$, then $e=m+1$ as $a\geq \frac{e-1-m}{2}$. Thus $\beta=2e$.
		
		(3) Suppose $S$ is a rational normal surface scroll $S(a,e-a)$ ($0\leq a\leq \frac{e}{2}$) and $C$ is nonhyperelliptic. 
		As before, write $C\equiv \alpha H+\beta F$ where $\alpha=0$ and $H:=0$ when $a=0$. 
		If $a\geq 1$, then $\alpha\geq 1$. 
		As $C$ is a nonhyperelliptic curve of genus $g=m\geq 3$, we have $\alpha \geq 3$.
		Hence $C$ has a trisecant line. 
		This is a contradiction as the homogeneous ideal of $C$ is generated by quadrics.
		Hence $a=0$ and we get $C\equiv (e+1+m)F$.
	\end{proof}
	
	Using \Cref{prop:unique embedding surface of minimal degree}, we can describe the divisor types of extremal varieties in the first hierarchy.
	
	\begin{thmalphabetmaintext}\label{thm:ACM with degree=e+u}
		Let $X \subseteq \p^r$ be a nondegenerate projective variety of dimension $n$, codimension $e$, and degree $d=e+1+m$ with $3 \leq m \leq e-1$. Then the following statements are equivalent:
		\begin{enumerate}
			\item[\textnormal{(i)}] $X$ is an ACM variety contained in \vmd[n+1].
			\item[\textnormal{(ii)}] One of the following holds:
			\begin{enumerate}
				\item[$(1)$] $e=4$, $d=8$, and $X$ lies on (a cone over) the Veronese surface as a divisor.
				\item[$(2)$] $X$ is contained in a unique $(n+1)$-fold rational normal scroll $Y = S(a_1 , a_2 , \ldots ,a_{n+1} )$ for some integers $0 \leq a_1 \leq a_2 \leq \ldots \leq a_{n+1}$ such that 
				\begin{enumerate}
					\item[$(\alpha)$] If $a_{n} \geq 1$ then the divisor class of $X$ is $2H+(1-e+m)F$  where $H$ and $F$ denote the hyperplane section and a ruling of $Y$, respectively.
					\item[$(\beta)$] If $a_{n} = 0$ then the divisor class of $X$ is $(e+1+m)F$  where $F \cong \p^n$ is the effective generator of the divisor class group of $Y$.
				\end{enumerate}
			\end{enumerate}
		\end{enumerate}
	\end{thmalphabetmaintext}
	
	\begin{proof}
		(i)$\Rightarrow$(ii)
		Let $Y$ be the \vmd[n+1] containing $X$.
		By taking general hyperplane sections, we get $C\subseteq S\subseteq\p^{e+1}$ where the curve $C$ and the surface $S$ are general sections of $X$ and $Y$, respectively.
		Then $C$ is a linearly normal curve of genus $g=m$.
		As $S$ also has the minimal degree, $S$ is the Veronese surface or a rational normal surface scroll $S(a,e-a)$. 
		When $S$ is the Veronese surface, $Y$ is not a rational normal scroll since if it is, then $S$ contains infinitely many projective lines. 
		Hence $Y$ is a cone over the Veronese surface.
		Then we get $e=4$, $d=8$.
		When $S$ is a rational normal surface scroll, $Y$ is also a rational normal scroll, and one may write $Y=S(a_1,\cdots,a_n,a_{n+1})$ $(0\leq a_1\leq\cdots\leq a_{n+1})$.
		If $a_n\geq 1$, then $X\equiv 2H+(1-e+m)F$ by \Cref{prop:unique embedding surface of minimal degree}.
		If $a_n=0$, then $a_{n+1}=e$ and $X\equiv (e+1+m)F$ by the same proposition.
		
		(ii)$\Rightarrow$(i) 
		Let $C\subseteq S\subseteq\p^{e+1}$ as before.
		If $S$ is the Veronese surface, then $C$ is ACM by \cite[Proposition 2.9]{MR1615938}.
		If $S$ is a singular rational normal surface scroll, then $C$ is ACM by \cite[Example 5.2]{MR1824228}.
		If $S$ is a smooth rational normal surface scroll, then $C$ is ACM as $C\equiv 2H+(1-e+m)F$ (cf. \cite[Theorem 5.10.(a)]{MR1615938}, \cite[\S 2]{MR2131135}, or \cite[Theorem 4.3]{MR3247023}). 
		Hence $C$ is ACM so that $X$ is also ACM.
		Then $\beta_{e-1,1}(X)=\beta_{e-1,1}(C)=e-1$ so that $X$ is contained in a \vmd[n+1].
	\end{proof}	
	
	\begin{remark}
		Let $X\subseteq\p^r$ be as in Theorem \ref{thm:ACM with degree=e+u}. Suppose that $X$ is a ACM variety contained in a \vmd[n+1]. Then $X$ is $3$-regular and the graded Betti numbers $\beta_{i,j}(X)$ of $X$ are completely determined as follows:
		\begin{equation*}
			\beta_{p,1}(X)=\begin{cases}
				p\binom{e+1}{p+1}-m\binom{e}{p-1} &\text{for $1\leq p\leq e+1-m$}\\
				p\binom{e}{p+1} &\text{for $e+1-m\leq p$}.
			\end{cases}
		\end{equation*} 
		and  
		\begin{equation*}
			\beta_{p,2}(X)=\begin{cases}
				0 & \text{for $1 \leq p \leq e-m$}\\
				\big(p-(e+1-m)\big)\binom{e}{e-p} & \text{for $e+1-m \leq p \leq e$}.
			\end{cases}
		\end{equation*}
		Indeed, $X$ is contained in a variety, denoted by $Z$, of minimal degree. Consider a general zero-dimensional section $\Gamma \subset\p^{e}$ of $X$ and a general curve section $D \subset\p^{e}$ of $Z$. In this situation, the set $\Gamma$ of $(e+1+m)$-points can be regarded as lying on a rational normal curve $D$. Note that $\beta_{i,j}(X) = \beta_{i,j}(\Gamma)$ for all $i,j$ since $X$ is ACM. We thus obtain the desired formula from
		\cite[Proposition 2.3]{MR1291122}. Note that $\beta_{p,1}$ appear as extremal cases in Theorem \ref{thm_upper_bound_general_maintext}.
	\end{remark}

	\bibliographystyle{amsalpha}
	\bibliography{ref.bib}
\end{document}